\input epsf
\documentstyle{amsppt}
\pagewidth{5.4truein}\hcorrection{0.55in}
\pageheight{8.4truein}\vcorrection{0.75in}
\TagsOnRight
\NoRunningHeads
\catcode`\@=11
\def\logo@{}
\footline={\ifnum\pageno>1 \hfil\folio\hfil\else\hfil\fi}
\topmatter
\title A dual of MacMahon's theorem on plane partitions
\endtitle
\endtopmatter
\document

\def\mysec#1{\bigskip\centerline{\bf #1}\message{ * }\nopagebreak\bigskip\par}

\def\myref#1{\item"{[{\bf #1}]}"}

\def\epf{\hfill{$\square$}\smallpagebreak}

\def\cite#1{\relaxnext@
  \def\nextiii@##1,##2\end@{[{\bf##1},\,##2]}%
  \in@,{#1}\ifin@\def\next{\nextiii@#1\end@}\else
  \def\next{[{\bf#1}]}\fi\next}
\def\proclaimheadfont@{\smc}

\define\M{\operatorname{M}}

\define\h{\operatorname{H}}
\def\H{H}

\define\twoline#1#2{\line{\hfill{\smc #1}\hfill{\smc #2}\hfill}}
\define\twolinetwo#1#2{\line{{\smc #1}\hfill{\smc #2}}}
\define\twolinethree#1#2{\line{\phantom{poco}{\smc #1}\hfill{\smc #2}\phantom{poco}}}

\def\mypic#1{\epsffile{#1}}

\bigskip
\bigskip
\twolinethree{{\smc Mihai Ciucu\footnote{Supported in part by NSF
  grant DMS-1101670.}}}{{\smc Christian
  Krattenthaler\footnote{Supported in part by the Austrian
Science Foundation FWF, grants Z130-N13 and S9607-N13,
the latter in the framework of the National Research Network
``Analytic Combinatorics and Probabilistic Number
Theory."}}}

\bigskip
\twolinethree{{\rm Indiana University}}{{\rm \rm Universit\"at Wien}}
\twolinethree{{\rm Department of Mathematics}}{{\rm Fakult\"at f\"ur Mathematik}}
\twolinethree{{\rm Bloomington, IN 47401, USA}}{{\rm Nordbergstra\ss e 15}}
\twolinethree{{\rm }}{{\rm A-1090 Wien, Austria}}



\bigskip
{\eightpoint {\smc Abstract.} {\rm A classical theorem of MacMahon states that the number of lozenge tilings of any centrally symmetric hexagon drawn on the triangular lattice is given by a beautifully simple product formula. In this paper we present a counterpart of this formula, corresponding to the {\it exterior} of a concave hexagon obtained by turning $120^\circ$ after drawing each side (MacMahon's hexagon is obtained by turning $60^\circ$ after each step).}}

\bigskip
\bigskip

\define\And{1}
\define\Andtwo{2}
\define\cekz{3}
\define\sc{4}
\define\ec{5}
\define\ov{6}
\define\ef{7}
\define\gd{8}
\define\anglepap{9}
\define\CLP{10}
\define\DT{11}
\define\Glaish{12}
\define\KaKZAC{13}
\define\Kuo{14}
\define\Kup{15}
\define\MacM{16}
\define\Sta{17}
\define\Ste{18}

\define\eaa{1.1}
\define\eab{1.2}
\define\eac{1.3}
\define\ead{1.4}
\define\eae{1.5}

\define\eba{2.1}
\define\ebb{2.2}

\define\eca{3.1}
\define\ecb{3.2}
\define\ecc{3.3}
\define\ecd{3.4}

\define\eda{4.1}
\define\edb{4.2}
\define\edc{4.3}
\define\edd{4.4}
\define\ede{4.5}
\define\edf{4.6}
\define\edg{4.7}
\define\edh{4.8}
\define\edi{4.9}
\define\edj{4.10}
\define\edk{4.11}
\define\edl{4.12}
\define\edm{4.13}
\define\edn{4.14}

\define\eea{5.1}
\define\eeb{5.2}
\define\eec{5.3}
\define\eed{5.4}
\define\eee{5.5}
\define\eef{5.6}
\define\eeg{5.7}
\define\eeh{5.8}

\define\taa{1.1}

\define\tba{2.1}
\define\tbb{2.2}

\define\tca{3.1}
\define\tcb{3.2}

\define\faa{1.1}
\define\fab{1.2}

\define\fba{2.1}
\define\fbb{2.2}
\define\fbc{2.3}
\define\fbd{2.4}

\define\fca{3.1}
\define\fcb{3.2}
\define\fcc{3.3}
\define\fcd{3.4}

\define\fdaa{4.1}
\define\fdab{4.2}
\define\fda{4.3}
\define\fdb{4.4}
\define\fdc{4.5}
\define\fdd{4.6}
\define\fde{4.7}
\define\fdf{4.8}

\define\fea{5.1}

\define\ffa{6.1}

\vskip-0.1in
\mysec{1. Introduction}

MacMahon's classical theorem \cite{\MacM} on the number of plane partitions that
fit in a given box 
(see \cite{\And}\cite{\Sta}\cite{\Andtwo}\cite{\Kup}\cite{\Ste}\cite{\KaKZAC} for more recent developments)
is equivalent to the fact that the number of lozenge\footnote{A {\it lozenge} is the union of two adjacent equilateral triangles
of side-length~1.} tilings
of a hexagon of side-lengths $a$, $b$, $c$, $a$, $b$, $c$ (in
cyclic order) on the triangular lattice is equal to  
$$
P(a,b,c):=
\frac{\h(a)\h(b)\h(c)\h(a+b+c)}
{\h(a+b)\h(a+c)\h(b+c)},
\tag\eaa
$$
where the hyperfactorials $\h(n)$ are defined by
$$
\h(n):=0!\,1!\cdots(n-1)!
\tag\eab
$$
(see Figure~{\faa} for an example). 

The hexagon is obtained by turning $60^\circ$ after drawing each side. If instead, after drawing each side one turns $120^\circ$ (the other natural amount of turning on the triangular lattice), one obtains a shape of the type illustrated in Figure~{\fab};  we will call such a shape a {\it shamrock}.

\topinsert
\twoline{\mypic{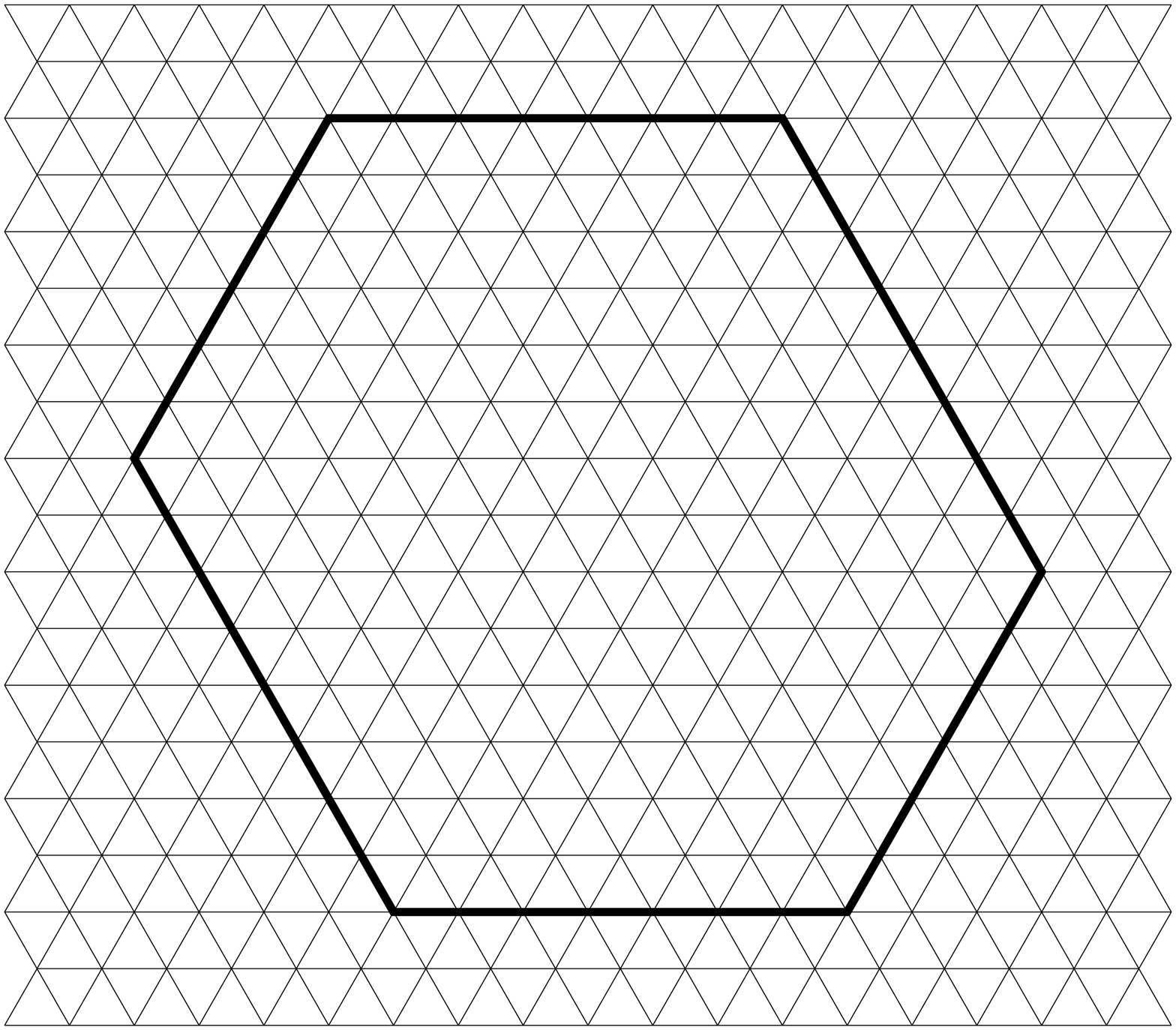}}{\mypic{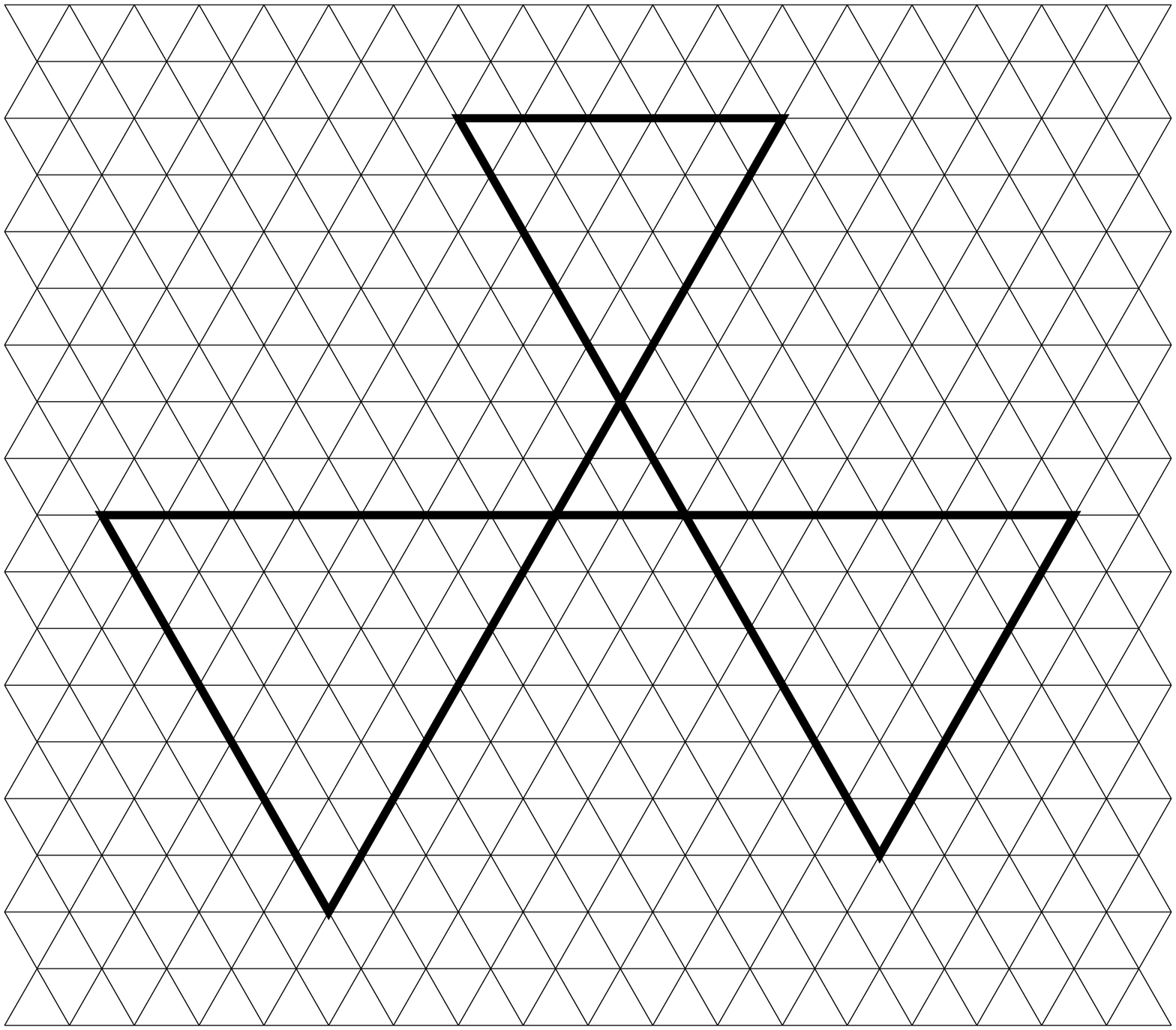}}
\twoline{Figure~{\faa}. {\rm Hexagon with $a=7$,}}{Figure~{\fab}. {\rm Shamrock with $m=2$,}}
\twoline{ {\rm \ \ \ \ \ \ \ \ \ \ \ \ \ $b=6$, $c=8$.}}{ {\rm \ \ \ \ \ \ \ \ \ \ \ \ \ \ \ \ \ \ \ \ \ \ \ $a=5$, $b=7$, $c=6$.}}
\endinsert

Our results concern the exterior of a shamrock (as its interior has no lozenge tilings). Let $S(a,b,c,m)$ be the shamrock whose central equilateral triangle has side-length $m$, while its top, bottom left and bottom right lobes are equilateral triangles of side-lengths $a$, $b$ and $c$, respectively; denote its exterior by $S^*(a,b,c,m)$.

We define the ratio of the number of tilings of the exteriors of the shamrocks $S(a,b,c,m)$ and $S(a+b+c,0,0,m)$ as follows. Let $H_N(a,b,c,m)$ be the hexagonal region of side-lengths alternating between $N+a+b+c$ and $N+a+b+c+m$ (the top side being $N+a+b+c$), and having the shamrock $S(a,b,c,m)$ removed from its center (to be precise, $H_N(a,b,c,m)$ is the region $SC_{N,N,N}(a,b,c,m)$ described in the next section). Then we define
$$
\frac{\M(S^*(a,b,c,m))}{\M(S^*(a+b+c,0,0,m))}
:=
\lim_{N\to\infty}\frac{\M(H_N(a,b,c,m))}{\M(H_N(a+b+c,0,0,m))},
\tag\eac
$$
where for a lattice region $R$, $\M(R)$ denotes the number of lozenge tilings of $R$.

The dual MacMahon theorem we obtain in this paper is the following. 

\proclaim{Theorem \taa} For any non-negative integers $a$, $b$, $c$ and $m$ we have
$$
\spreadlines{3\jot}
\align
\frac{\M(S^*(a,b,c,m))}{\M(S^*(a+b+c,0,0,m))}
&=
\frac{\h(a)\h(b)\h(c)\h(a+b+c+m)\h(m)^2}{\h(a+m)\h(b+m)\h(c+m)\h(a+b+c)}
\\
&=P(a,b,m)\,P(a+b,c,m).
\tag\ead
\endalign
$$
In the special case when $m=a+b+c$, this can be written simply as
$$
\frac{\M(S^*(a,b,c,a+b+c))}{\M(S^*(a+b+c,0,0,a+b+c))}
=P(a,b,c)\,P(a+b,b+c,c+a).
\tag\eae
$$

\endproclaim

Formula (\eae) is illustrated geometrically in Figure~{\fea} (see
Section~5). Theorem~{\taa} will follow as a consequence of a more
general result (see Theorems~{\tba} and {\tbb}), which we describe 
in the next section. The purpose of Section~3 is to establish
enumeration formulas for so-called ``magnet bar regions,"
which will provide the base cases for the inductive proof of
Theorems~{\tba} and {\tbb} given in Section~4. Finally, in
Section~5, we show how Theorem~{\taa} follows from Theorems~{\tba} and
{\tbb}. There, we also explain how the results of this paper
relate to work of the first author on correlation of holes ``in a sea
of dimers," and we comment on the physical interpretation of our
results. We close the paper by briefly highlighting the 
achievements of this paper when compared to the earlier paper
\cite{\cekz}, and what perspectives they offer.

\mysec{2. Precise statement of results}

The family of regions whose lozenge tilings we enumerate in this paper is a generalization of the following family, introduced in \cite{\cekz}. Consider hexagons of sides $x$, $y+m$, $z$, $x+m$, $y$, $z+m$ (in clockwise order, starting from top) with an equilateral triangle of side $m$ removed from its center (see Figures~{\fba} and {\fbb} for examples). This triangle is called the {\it core}, and the leftover region, denoted $C_{x,y,z}(m)$, a {\it cored hexagon}.

To define $C_{x,y,z}(m)$ precisely, we need to specify what position of the core is the``central'' one. Let $s$ be a side of the core, and let $u$ and $v$ be the sides of the hexagon parallel to it. The most natural definition would require that the distance between $s$ and $u$ is the same as the distance between $v$ and the vertex of the core opposite $s$, for all three choices of $s$.

However, since the sides of the core have to be along lines of the
underlying triangular lattice, it is easy to see that this can be
achieved only if $x$, $y$ and $z$ have the same parity (Figure~{\fba}
illustrates such a case); in that case, we define this to be the
position of the core. On the other hand, if for instance $x$ has
parity different from that of $y$ and $z$, the triangle satisfying the
above requirements would only have one side along a lattice line,
while each of the remaining two extends midway between two consecutive
lattice lines (this can be seen from Figure~{\fbb}). To resolve this,
we translate this central triangle half a unit towards the side of the
hexagon of length $y$, in a direction parallel to the side of length
$x$, and define this to be the position of the core in this case
(see Figure~{\fbb}).

\topinsert
\twoline{\mypic{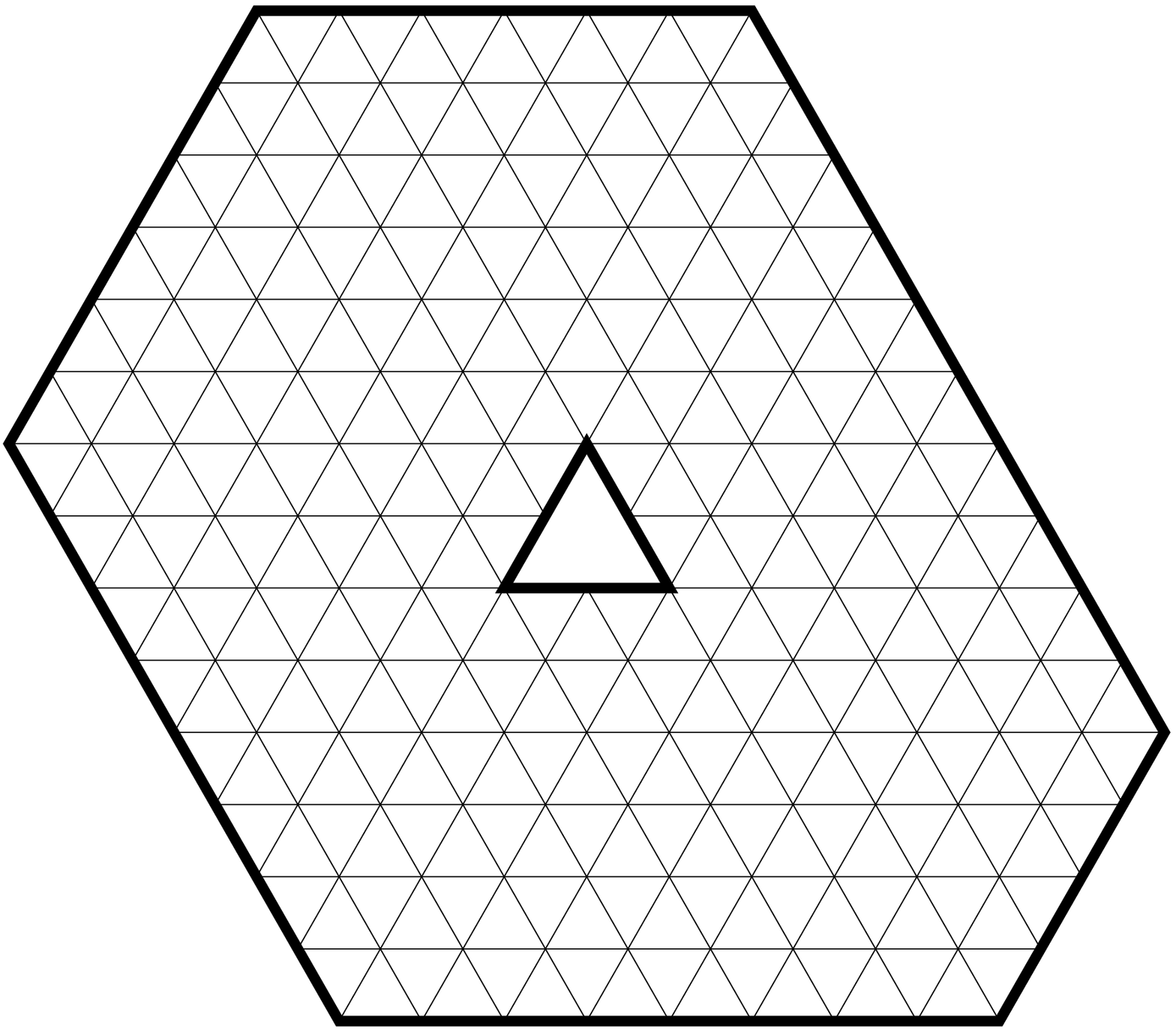}}{\mypic{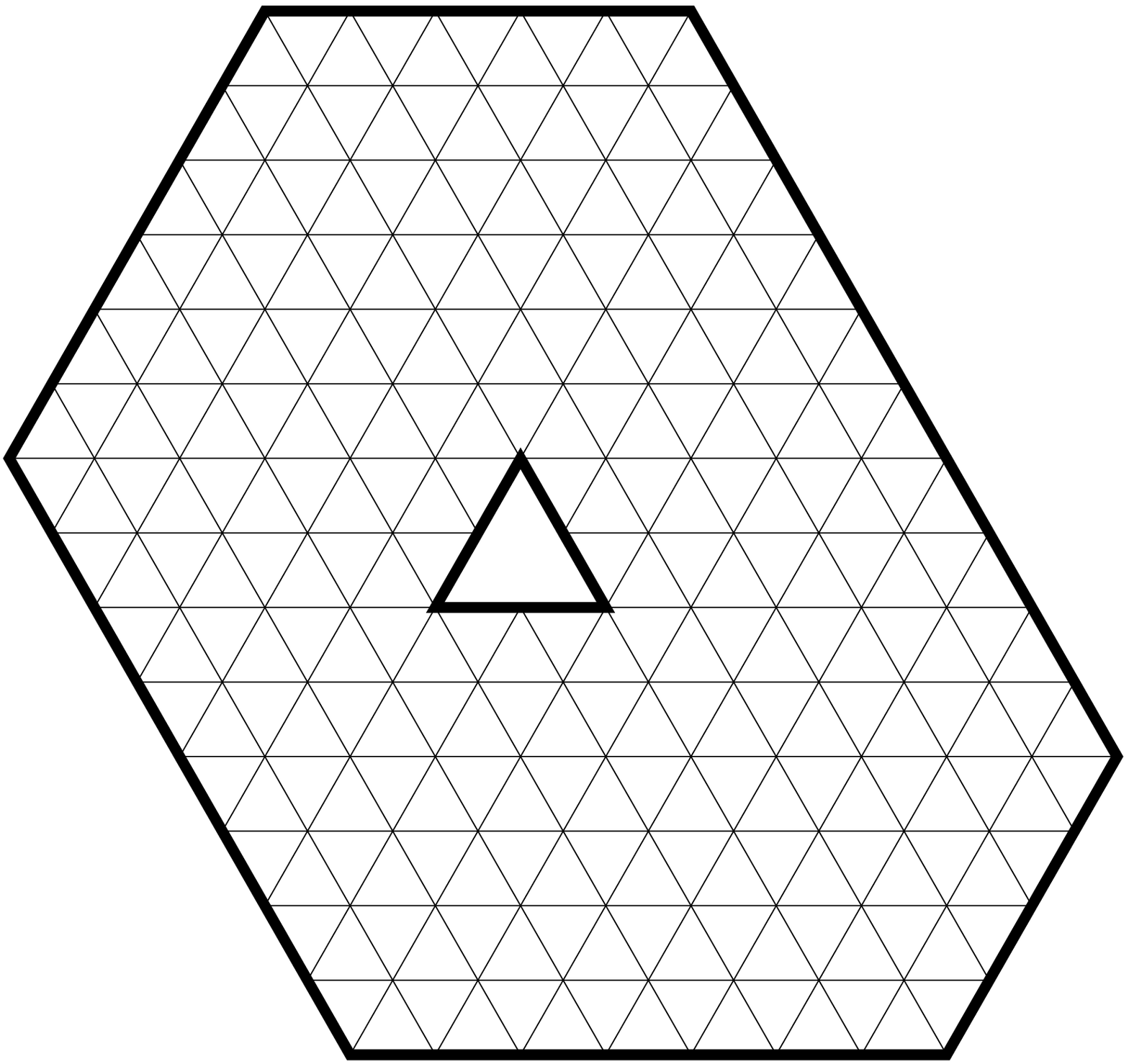}}
\twoline{Figure~{\fba}.{\rm \ $C_{6,8,4}(2)$.}}
{Figure~{\fbb}.{\rm \ $C_{5,8,4}(2)$.}}
\endinsert

The family of regions we will be concerned with in this paper is a generalization of cored hexagons, corresponding to the case when the core is not just a triangle, but a shamrock.

Given non-negative integers $a$, $b$ and $c$, construct our region as follows. Start with the cored hexagon $C_{x,y,z}(m)$, and ``push out'' the six lattice lines containing its sides as follows: the top, southwestern and southeastern sides $a$, $b$ and $c$ units, respectively, and the bottom, northeastern and northwestern sides $b+c$, $a+c$ and $a+b$ units, respectively. Enlarge also the core by adding to it downpointing equilateral triangles of sides $a$, $b$ and $c$ that touch the original core at its top, left and right vertex, respectively. The hexagon formed by the pushed out edges, with the enlarged, shamrock-shaped core taken out of it, is called an {\it S-cored hexagon}, and is denoted by $SC_{x,y,z}(a,b,c,m)$ (see Figures~{\fbc} and {\fbd} for two examples; the dotted lines indicate the underlying cored hexagons).

\topinsert
\twoline{\mypic{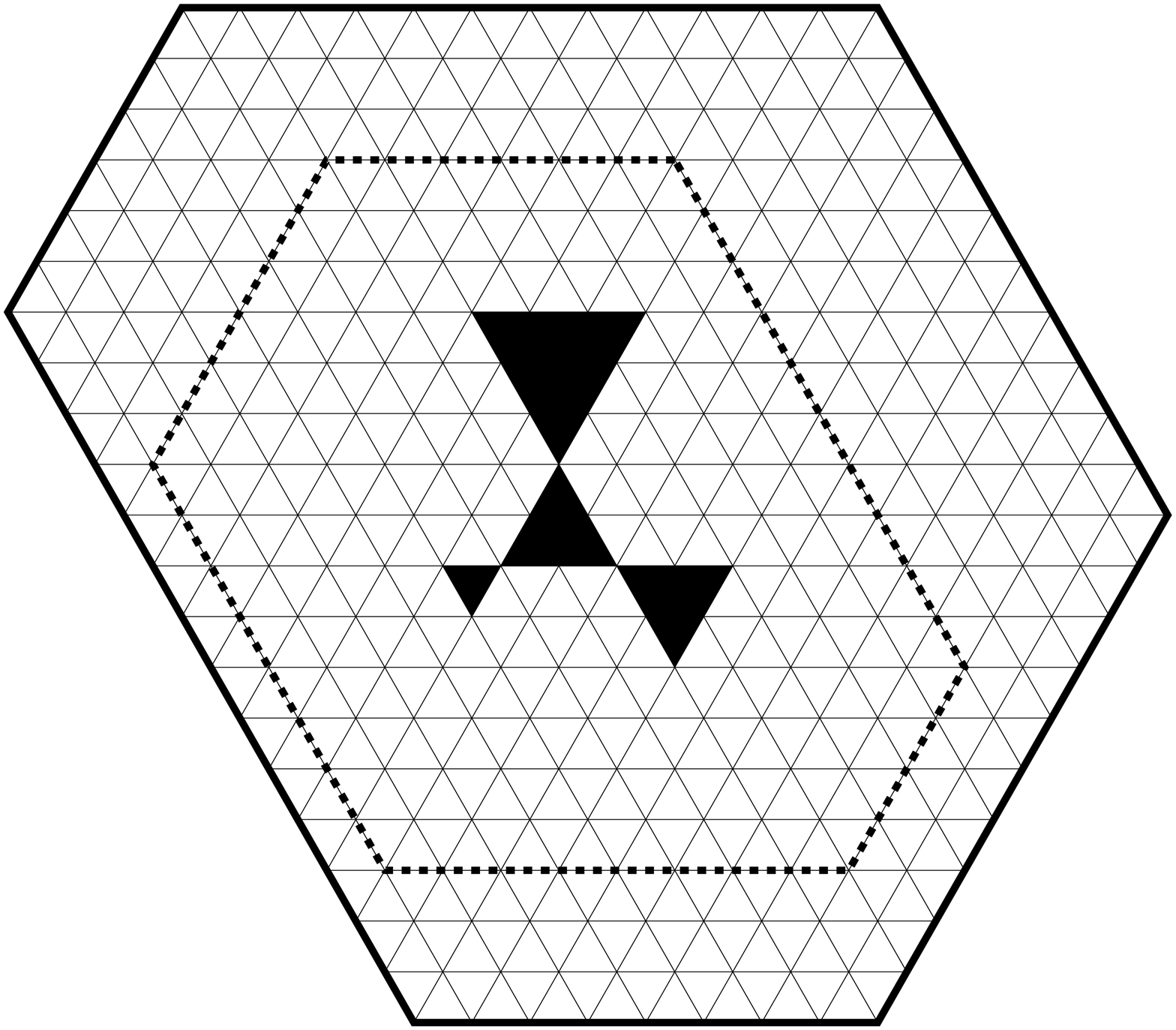}}{\mypic{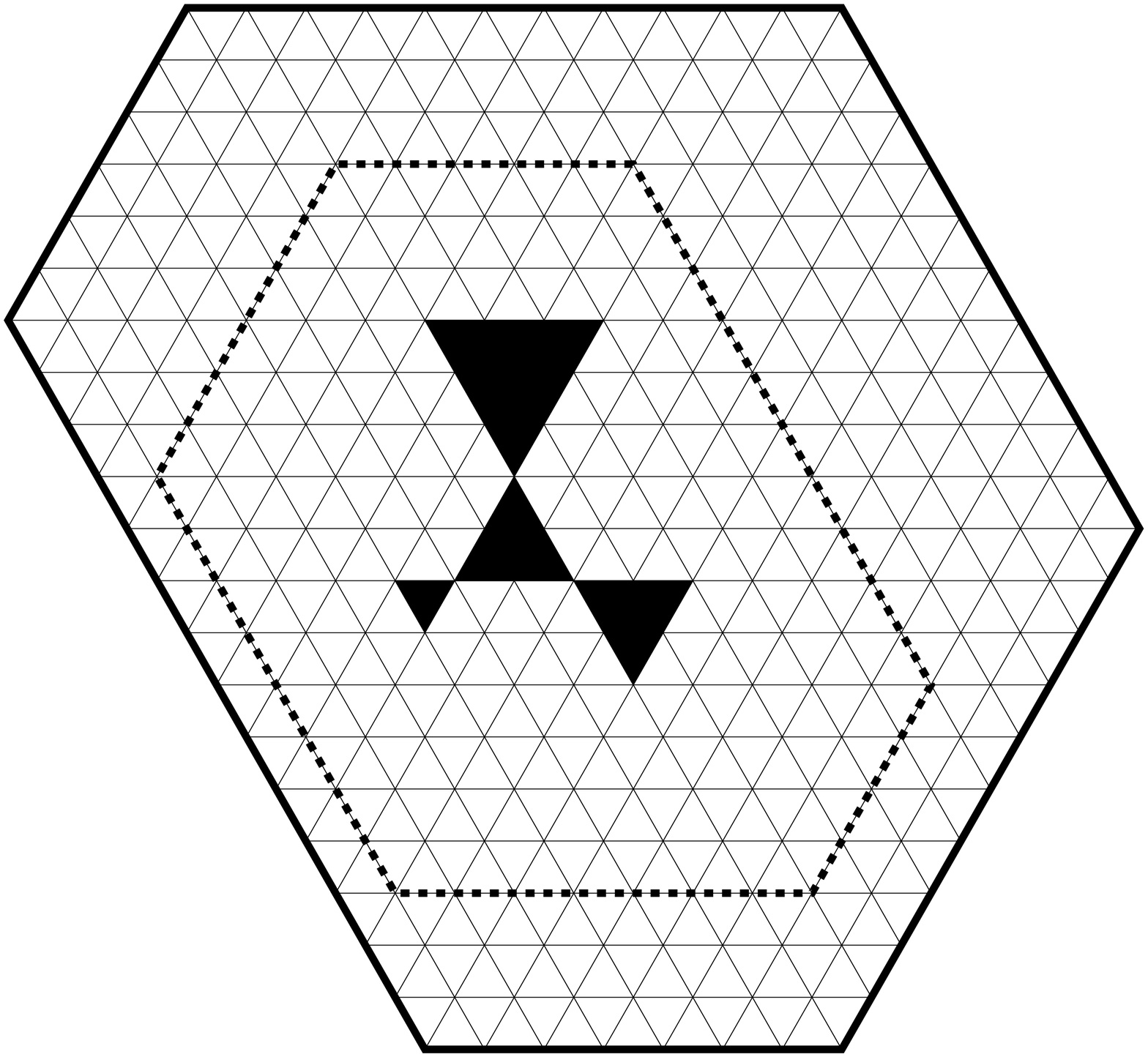}}
\twoline{Figure~{\fba}.{\rm \ $SC_{6,8,4}(3,1,2,2)$.}}
{Figure~{\fbb}.{\rm \ $SC_{5,8,4}(3,1,2,2)$.}}
\endinsert

The main result of this paper --- from which Theorem~{\taa} will follow easily --- is the exact enumeration of the lozenge tilings of $S$-cored hexagons. This is contained in the following two results.
In order to formulate them, it is convenient to extend the
definition of hyperfactorials to half-integers (i.e.,
odd integers divided by~2):
$$\h(n):=\cases
\prod _{k=0} ^{n-1}{\Gamma(k+1)},\quad &\text {for $n$ a positive integer,}\\
\prod _{k=0} ^{n-\frac {1} {2}}{\Gamma(k+\frac {1} {2})}, \quad &\text 
{for $n$ a positive half-integer},
\endcases
$$ 
where $\Gamma$ denotes the classical gamma function.

\proclaim{Theorem \tba}  Let $x$, $y$, $z$, $a$, $b$, $c$ and $m$ be nonnegative integers. If $x$, $y$ and $z$ have the same parity, we have
$$
\spreadlines{0\jot}
\align
&
\M(SC_{x,y,z}(a,b,c,m))=
\frac {
\h(m)^3\,
\h(a)\,
\h(b)\,
\h(c)
} 
{
\h(m+a)\,
\h(m+b)\,
\h(m+c)
}
\\
&
\times
\frac 
{
\h(\frac{x+y}2+m+a+b)\,
\h(\frac{x+z}2+m+a+c)\,
\h(\frac{y+z}2+m+b+c)
} 
{
\h(\frac{x+y}2+m+c)\,
\h(\frac{x+z}2+m+b)\,
\h(\frac{y+z}2+m+a)
}
\\
&
\times
\frac 
{
\h(\frac{x+y}2+c)\,
\h(\frac{x+z}2+b)\,
\h(\frac{y+z}2+a)
}
{
\h(\frac{x+y}2+a+b)\,
\h(\frac{x+z}2+a+c)\,
\h(\frac{y+z}2+b+c)
} 
\\
&
\times
\frac {\h(x + m+a+b+c)\h(y + m+a+b+c)
} 
{\h(x + y + m+a+b+c)\h(x + z + m+a+b+c)
}\\
&
\times
\frac {\h(z + m+a+b+c)\h(x + y + z + m+a+b+c)} 
{\h(y + z + m+a+b+c)}\\
&
\times
\frac {
\h(\lceil{\frac{x + y + z}2}\rceil+m+a+b+c )
\h(\lfloor{\frac{x + y + z}2}\rfloor+m+a+b+c )
} 
{
\h(\frac{x + y}2+m+a+b+c)     
\h(\frac{x + z}2+m+a+b+c)
\h(\frac{y + z}2+m+a+b+c)
}
\\
&
\times\frac {
\h(\lceil{\frac{x}2}\rceil
\h(\lfloor{\frac{x}2}\rfloor)
\h(\lceil{\frac{y}2}\rceil)
} 
{
\h(\lceil{\frac{x}2}\rceil+{\frac{m+a+b+c}{2}})
\h(\lfloor{\frac{x}2}\rfloor+{\frac{m+a+b+c}{2}})
     \h(\lceil{\frac{y}2}\rceil+{\frac{m+a+b+c}{2}})
     }
\\
&
\times\frac {
\h(\lfloor{\frac{y}2}\rfloor)
     \h(\lceil{\frac{z}2}\rceil)
     \h(\lfloor{\frac{z}2}\rfloor)
} 
{
     \h(\lfloor{\frac{y}2}\rfloor+{\frac{m+a+b+c}{2}})
     \h(\lceil{\frac{z}2}\rceil+{\frac{m+a+b+c}{2}})
     \h(\lfloor{\frac{z}2}\rfloor+{\frac{m+a+b+c}{2}})
}
\\
&
\times
\frac {\h(\frac{m+a+b+c}{2})^2 \,
\h(\frac{x+y}2+{\frac{m+a+b+c}{2}})^2\,
\h(\frac{x+z}2+{\frac{m+a+b+c}{2}})^2\,
\h(\frac{y+z}2+{\frac{m+a+b+c}{2}})^2
} 
{
\h(\lceil{\frac{x + y + z}2}\rceil+{\frac{m+a+b+c}{2}})
\h(\lfloor{\frac{x + y + z}2}\rfloor+{\frac{m+a+b+c}{2}})
\h(\frac{x + y}2)\h(\frac{x + z}2)\h(\frac{y + z}2)
               }.
\tag\eba
\endalign
$$

\endproclaim

\proclaim{Theorem \tbb} Let $x$, $y$, $z$, $a$, $b$, $c$ and $m$ be nonnegative integers. If $x$ has parity different from the parity of $y$ and $z$, we have
$$
\align
&
\M(SC_{x,y,z}(a,b,c,m))=
\frac {
\h(m)^3\,
\h(a)\,
\h(b)\,
\h(c)
} 
{
\h(m+a)\,
\h(m+b)\,
\h(m+c)
}
\\
&
\times
\frac 
{
\h(\lfloor{\frac{x+y}2}\rfloor+m+a+b)\,
\h(\lceil{\frac{x+z}2}\rceil+m+a+c)\,
\h(\frac{y+z}2+m+b+c)
} 
{
\h(\lceil{\frac{x+y}2}\rceil+m+c)\,
\h(\lfloor{\frac{x+z}2}\rfloor+m+b)\,
\h(\frac{y+z}2+m+a)
}
\\
&
\times
\frac 
{
\h(\lceil{\frac{x+y}2}\rceil+c)\,
\h(\lfloor{\frac{x+z}2}\rfloor+b)\,
\h(\frac{y+z}2+a)
}
{
\h(\lfloor{\frac{x+y}2}\rfloor+a+b)\,
\h(\lceil{\frac{x+z}2}\rceil+a+c)\,
\h(\frac{y+z}2+b+c)
} 
\\
&
\times
\frac {\h(x + m+a+b+c)\h(y + m+a+b+c)
} 
{\h(x + y + m+a+b+c)\h(x + z + m+a+b+c)
}\\
&
\times
\frac {\h(z + m+a+b+c)\h(x + y + z + m+a+b+c)} 
{\h(y + z + m+a+b+c)}\\
&
\times
\frac {
\h(\lceil{\frac{x + y + z}2}\rceil+m+a+b+c )
\h(\lfloor{\frac{x + y + z}2}\rfloor+m+a+b+c )
} 
{
\h(\lfloor{\frac{x + y}2}\rfloor+m+a+b+c)     
\h(\lceil{\frac{x + z}2}\rceil+m+a+b+c)
\h(\frac{y + z}2+m+a+b+c)
}
\\
&
\times\frac {
\h(\lceil{\frac{x}2}\rceil)
\h(\lfloor{\frac{x}2}\rfloor)
\h(\lceil{\frac{y}2}\rceil)
} 
{
\h(\lceil{\frac{x}2}\rceil+{\frac{m+a+b+c}{2}})
\h(\lfloor{\frac{x}2}\rfloor+{\frac{m+a+b+c}{2}})
     \h(\lceil{\frac{y}2}+{\frac{m+a+b+c}{2}})
     }
\\
&
\times\frac {
\h(\lfloor{\frac{y}2}\rfloor)
     \h(\lceil{\frac{z}2}\rceil)
     \h(\lfloor{\frac{z}2}\rfloor)
} 
{
     \h(\lfloor{\frac{y}2}\rfloor+{\frac{m+a+b+c}{2}})
     \h(\lceil{\frac{z}2}\rceil+{\frac{m+a+b+c}{2}})
     \h(\lfloor{\frac{z}2}\rfloor+{\frac{m+a+b+c}{2}})
}
\\
&
\times
\frac {\h(\frac{m+a+b+c}{2})^2 \,
\h(\lceil{\frac{x+y}2}\rceil+{\frac{m+a+b+c}{2}})
\h(\lfloor{\frac{x+y}2}\rfloor+{\frac{m+a+b+c}{2}})
} 
{
\h(\lceil{\frac{x + y + z}2}\rceil+{\frac{m+a+b+c}{2}})
\h(\lfloor{\frac{x + y + z}2}\rfloor+{\frac{m+a+b+c}{2}})
\h(\lceil{\frac{x + y}2}\rceil)
               }\\
&
\times
\frac {
\h(\lceil{\frac{x+z}2}\rceil+{\frac{m+a+b+c}{2}})
\h(\lfloor{\frac{x+z}2}\rfloor+{\frac{m+a+b+c}{2}})
\h(\frac{y+z}2+{\frac{m+a+b+c}{2}})^2
} 
{
\h(\lfloor{\frac{x + z}2}\rfloor)\h(\frac{y + z}2)
               }.
\tag\ebb
\endalign
$$

\endproclaim

This represents a common generalization of MacMahon's formula~(\eaa) and Theorem~{\taa}. It shows that, if one regards the hexagon as being the right outer boundary  to consider on the triangular lattice for the corresponding region to have a number of lozenge tilings given by a simple product formula, then a good inner boundary is the shamrock.

Note also that Theorems~{\tba} and {\tbb} generalize the main results of \cite{\cekz}, by introducing three new parameters to the geometry of the core (the sizes of the three lobes of the shamrock). This results in a new, four parameter generalization of MacMahon's theorem~(\eaa).

\mysec{3. Magnet bar regions}

Our proof of Theorems~{\tba} and {\tbb} will use the exact enumeration of the following regions.

\topinsert
\centerline{\mypic{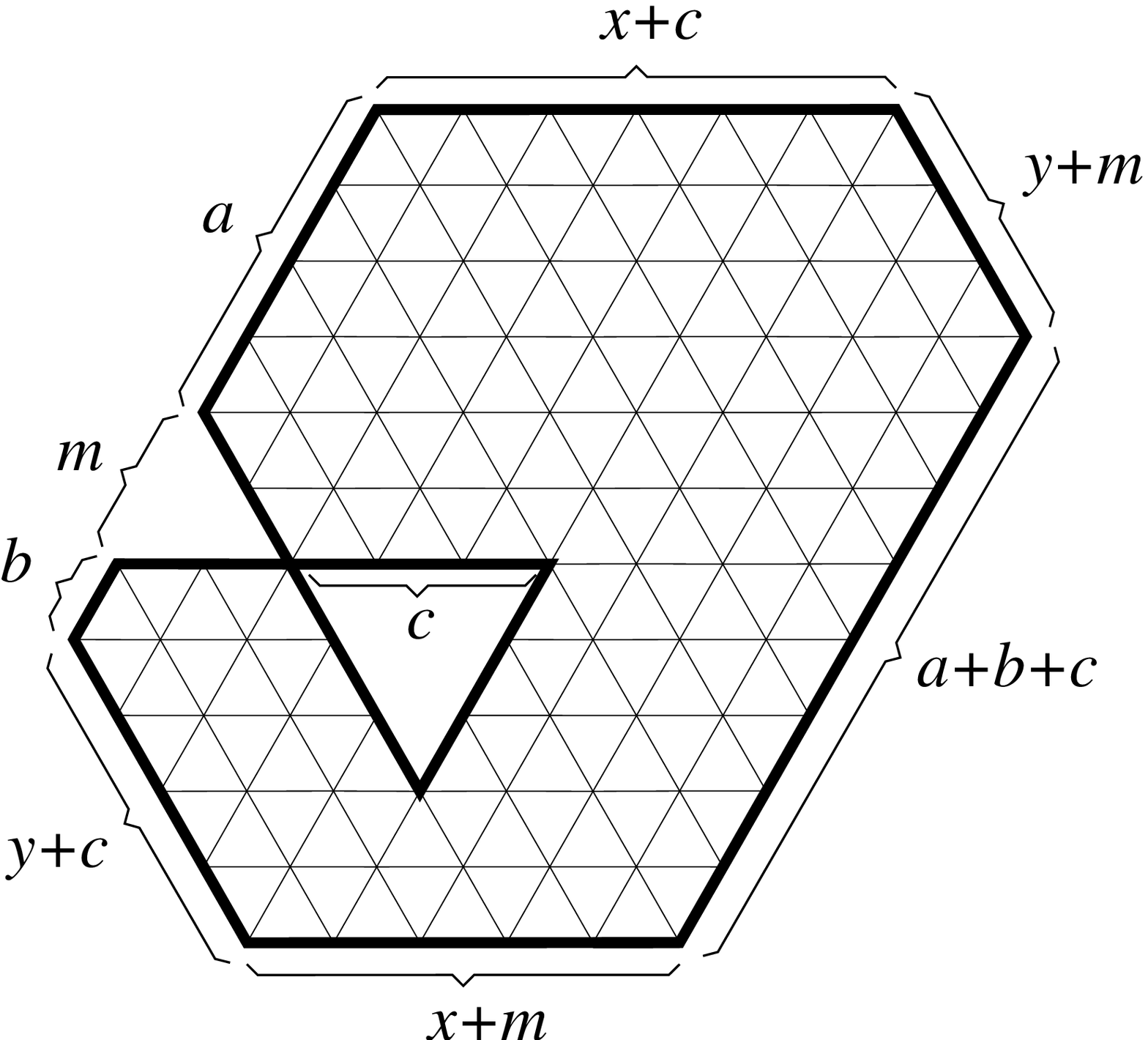}}
\medskip
\centerline{{\smc Figure~{\fca}{\rm . The magnet bar region $B_{3,1}(4,1,3,2)$.}}}
\endinsert

Let $x$, $y$, $a$, $b$, $c$ and $m$ be nonnegative integers. Consider the hexagonal region with two removed equilateral triangles illustrated by Figure~{\fca}, where the top side lengths of the hexagon are, clockwise from top, $x+c$, $y+m$, $a+b+c$, $x+m$, $y+c$, $a+b+m$, and the sides of the two removed triangles are $m$ (for the triangle touching the northwestern hexagon side; note that the lengths of the portions of that side above and below this triangle are $a$ and $b$, respectively) and $c$ (for the central triangle). We call such a region a {\it magnet bar region} and denote it by $B_{x,y}(a,b,c,m)$.

The main result of this section is the following.

\proclaim{Theorem \tca} For any nonnegative integers $x$, $y$, $a$, $b$, $c$ and $m$, we have
$$
\spreadlines{3\jot}
\align
\M(B_{x,y}(a,b,c,m))
&
=
\frac{\h(m)^2\,\h(a)\,\h(b)\,\h(c)\,\h(m+a+b+c)}
{\h(m+a)\,\h(m+b)\,\h(m+c)}
\\
&
\times
\frac{\h(x+m+a+c)\,\h(y+m+b+c)}{\h(x+y+m+c)}
\frac{\h(x+y+c)}{\h(x+a+c)\,\h(y+b+c)}
\\&
\times
\frac{\h(x+y+m+a+b+c)}{\h(x+m+a+b+c)\,\h(y+m+a+b+c)}
\frac{\h(x)\,\h(y)}{\h(x+y)}.
\tag\eca
\endalign
$$

\endproclaim

The proof of the above result (as well as the proof of Theorems~{\tba}
and {\tbb}) is based on Kuo's powerful graphical condensation method
(see \cite{\Kuo}). 
For ease of reference, we state below the particular instance of Kuo's general results that we need for our proofs (which is Theorem~2.1 in \cite{\Kuo}).

\proclaim{Theorem {\tcb} \smc(Kuo)} Let $G=(V_1,V_2,E)$ be a plane
bipartite graph in which $|V_1|=|V_2|$. Let vertices 
$\alpha$, $\beta$, $\gamma$ and $\delta$ appear cyclically on a face of $G$. If $\alpha,\gamma\in V_1$ and $\beta,\delta\in V_2$, then
$$
\M(G)\M(G-\{\alpha,\beta,\gamma.\delta\})=
\M(G-\{\alpha,\beta\})\M(G-\{\gamma,\delta\})+
\M(G-\{\alpha,\delta\})\M(G-\{\beta,\gamma\}).
\tag\ecb
$$

\endproclaim

{\it Proof of Theorem {\tca}.} We prove (\eca) by induction on $x+y+b$. Our base cases will be the situations when $x=0$, $y=0$ or $b=0$. If $b=0$, note that the lozenge tiling is forced in an $m\times(y+c)$ rectangle along the southwestern edge of the magnet bar region (Figure~{\fcb} illustrates this for $x=3$, $y=1$, $a=4$, $c=3$, $m=2$). Upon removing these forced lozenges, the leftover region is one whose tilings are enumerated by a formula due to Cohn, Larsen and Propp (see Proposition~2.1 in \cite{\CLP}). Indeed, regard the leftover region --- shown in bold contour in Figure~{\fcb} --- as being obtained from a hexagon by cutting out from it a triangle of side $c$ resting on its southwestern edge. It is readily seen that this region has the same number of lozenge tilings as the region obtained from the hexagon by cutting out $c$ consecutive unit triangles resting on the southwestern edge (this is due to forced lozenges in any tiling of the latter). However, regions of this type have their number of lozenge tilings given by \cite{\CLP, Proposition~2.1}. It is routine to check that the resulting formula agrees with the $b=0$ specialization of the expression on the right-hand side of (\eca).

\topinsert
\twoline{\mypic{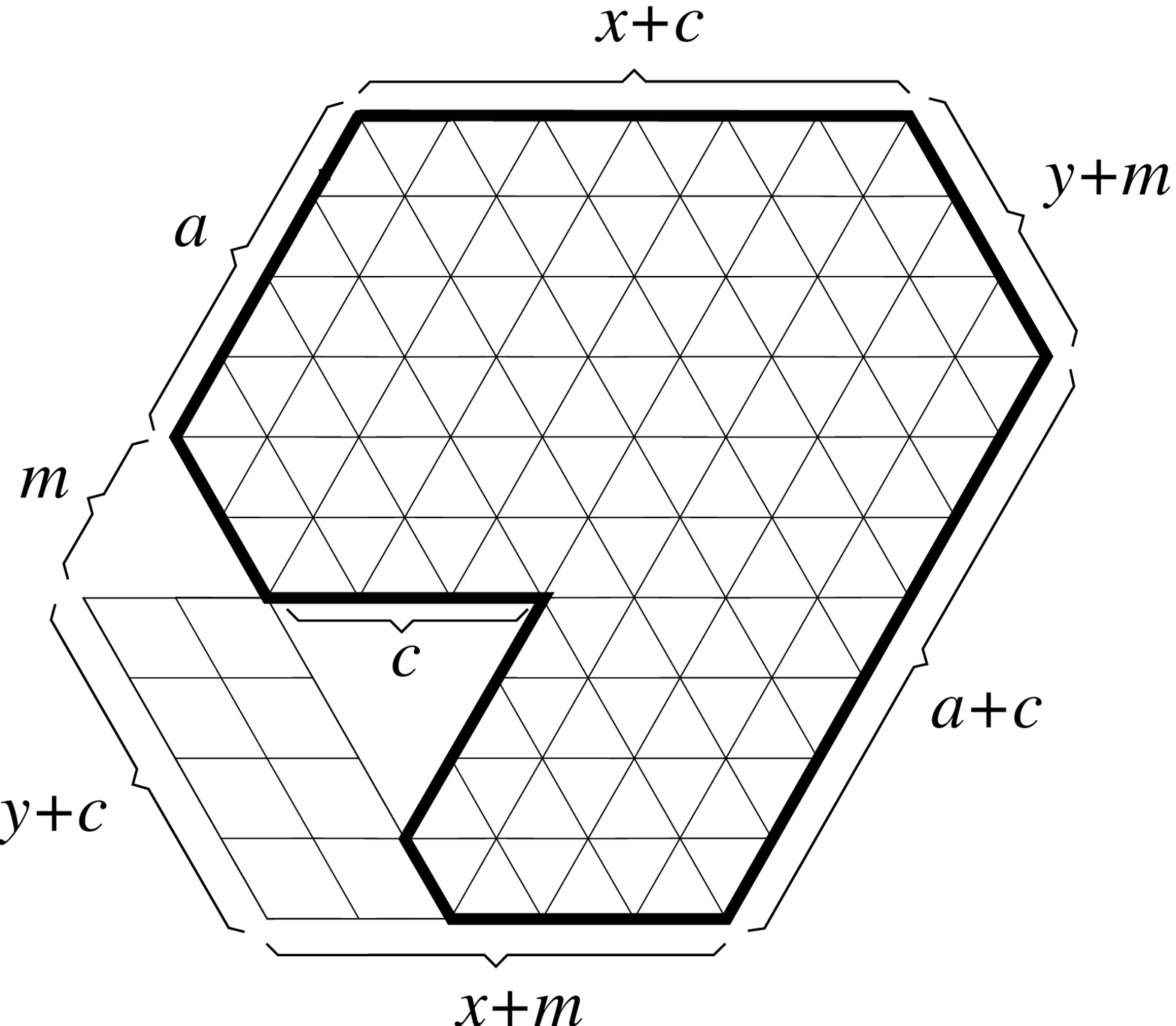}}{\mypic{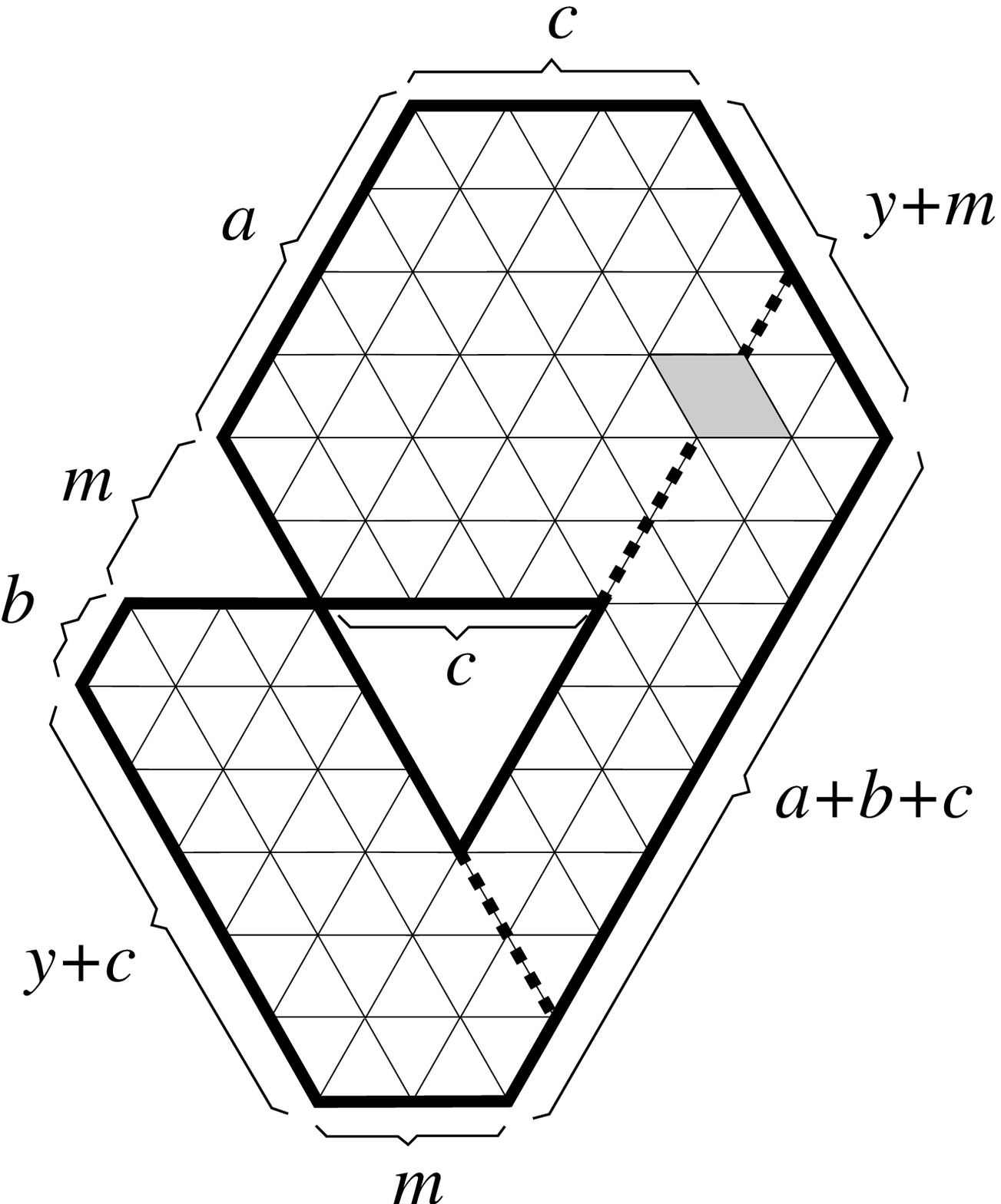}}
\medskip
\twoline{Figure~{\fcb}.{\rm \ The case $b=0$.}}
{Figure~{\fcc}.{\rm \ The case $x=0$.}}
\endinsert

If $x=0$, it turns out that tiling the magnet bar region with lozenges
is equivalent to tiling two disjoint, hexagon shaped subregions (this
is illustrated in Figure~{\fcc}; the hexagonal subregions are
determined by the two dashed lines;
that one of them is interrupted by
a shaded lozenge should be ignored at this point). 

To see this, consider for instance the top hexagon 
$\H_1$, obtained by extending upward the southeastern side of the triangle of side $c$ (shown in Figure~{\fcc} by a dashed line). Let $T$ be a lozenge tiling of the magnet bar region. We claim that there is no lozenge in $T$ that straddles the southeastern side of $\H_1$. Indeed, suppose towards a contradiction that $T$ contained such a lozenge $L$ (pictured in gray in Figure~{\fcc}). Following upward from $L$ the lozenge in $T$ above $L$, then the lozenge in $T$ above that, and so on, one obtains a path of lozenges that eventually must end on the upper side of the magnet bar region. However, the top side of the cut out triangle of side $c$ generates in a similar fashion $c$ more paths of rhombi, which must also end at the top side of the magnet bar region. Since the latter side has length $c$ and the above $c+1$ paths of rhombi are necessarily nonintersecting (as they are all part of the tiling $T$), this provides the contradiction that proves our claim. 

Thus the top dashed line cannot be crossed by any lozenge. It follows that the rhombus cut out from the eastern corner by the two dashed lines is forced to be part of each tiling, and can therefore be removed without affecting the number of tilings of the magnet bar region. The leftover region is the union of two disjoint hexagons. Using (\eaa) it follows that for $x=0$ we have
$$
\M(B_{0,y}(a,b,c,m))=P(a,c,m)\,P(b,y+c,m).\tag\ecc
$$
It is apparent from (\eaa) that the right-hand side above agrees with
the $x=0$ specialization of the right-hand side of (\eca). The case 
$y=0$ follows similarly. This completes the verification of the base cases of our induction.

The induction step is based on a convenient application of Kuo's graphical condensation stated in Theorem~{\tcb}.

\topinsert

\twoline{\mypic{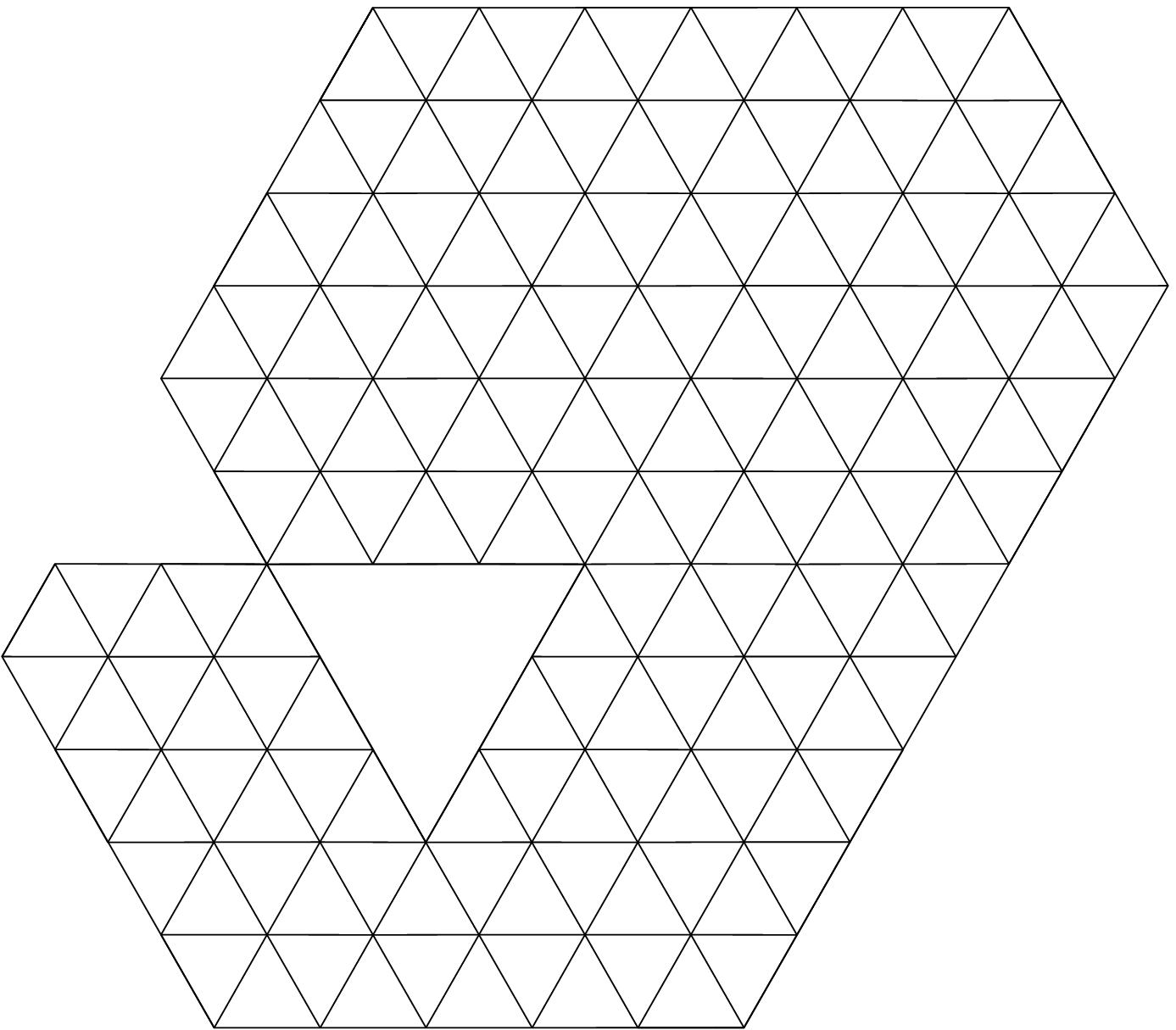}}{\mypic{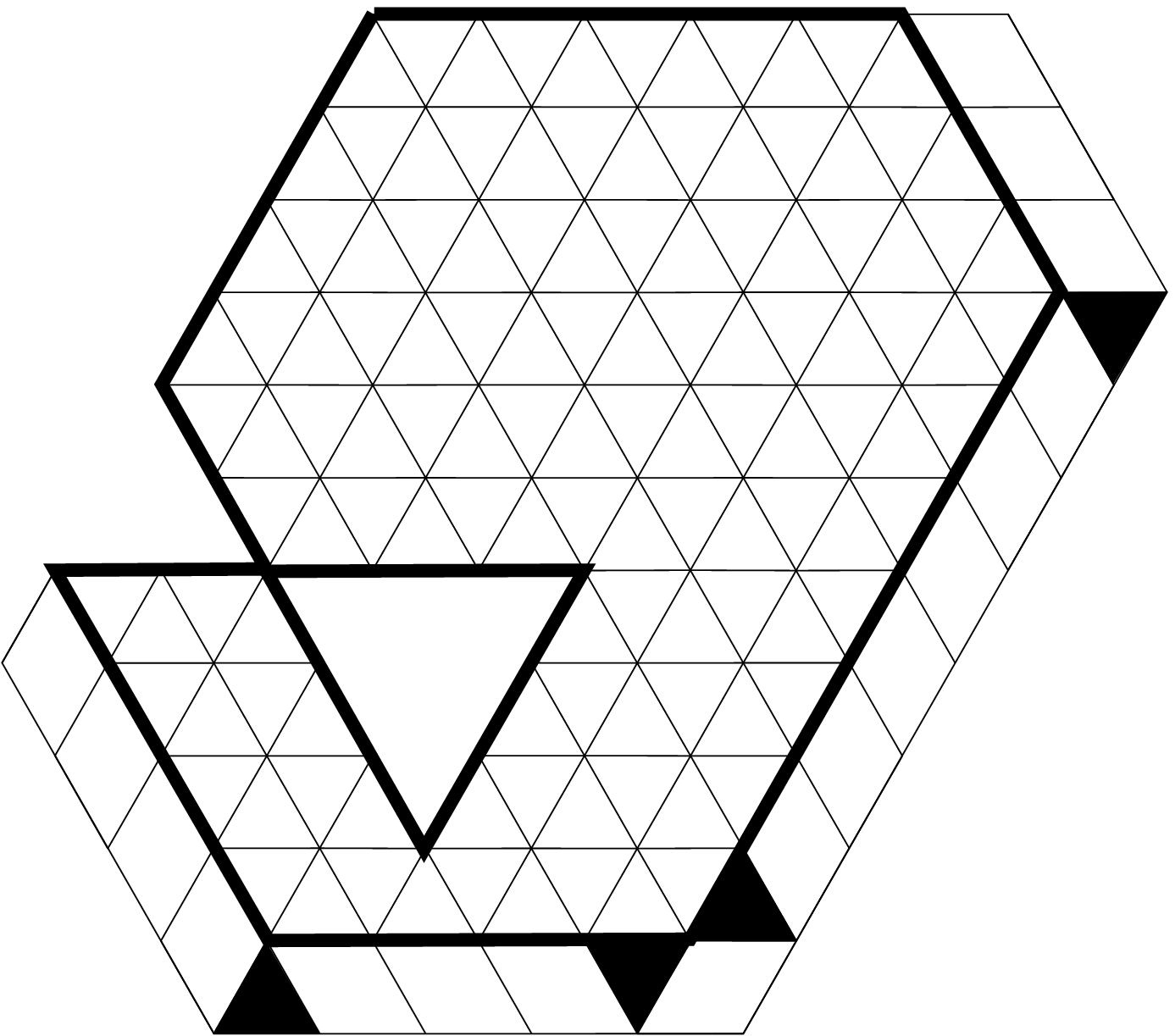}}
\bigskip

\twoline{\mypic{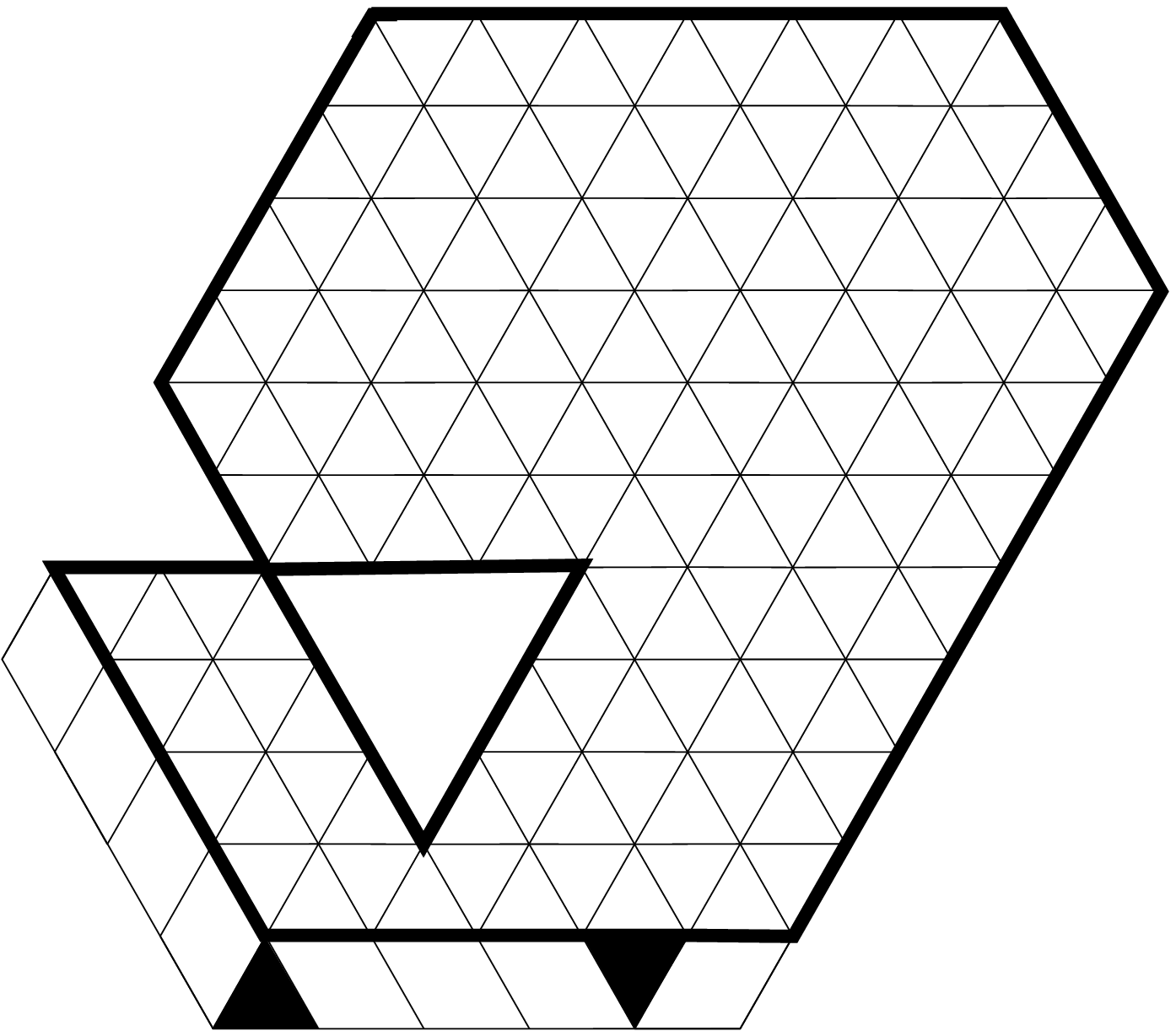}}{\mypic{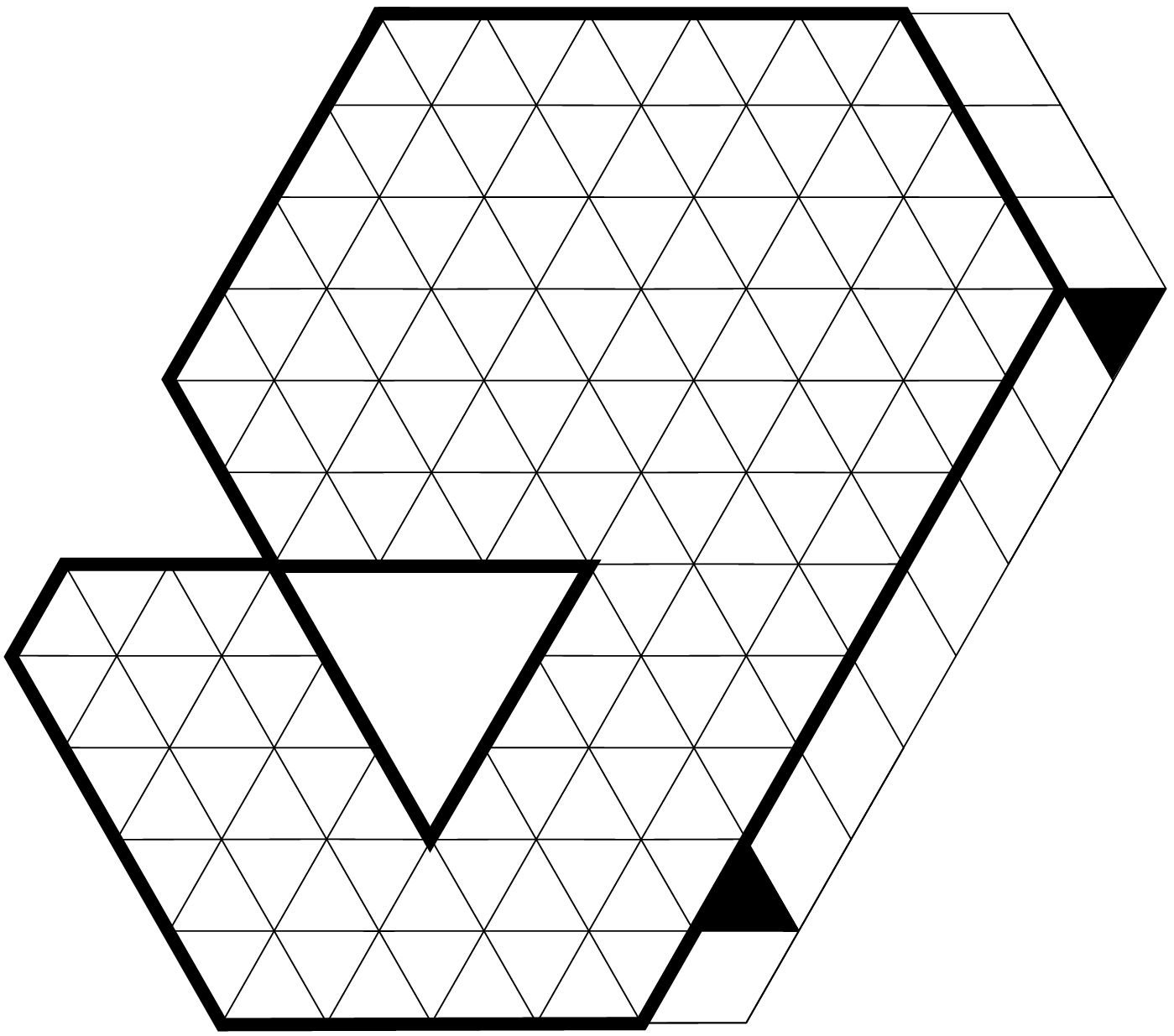}}
\bigskip

\twoline{\mypic{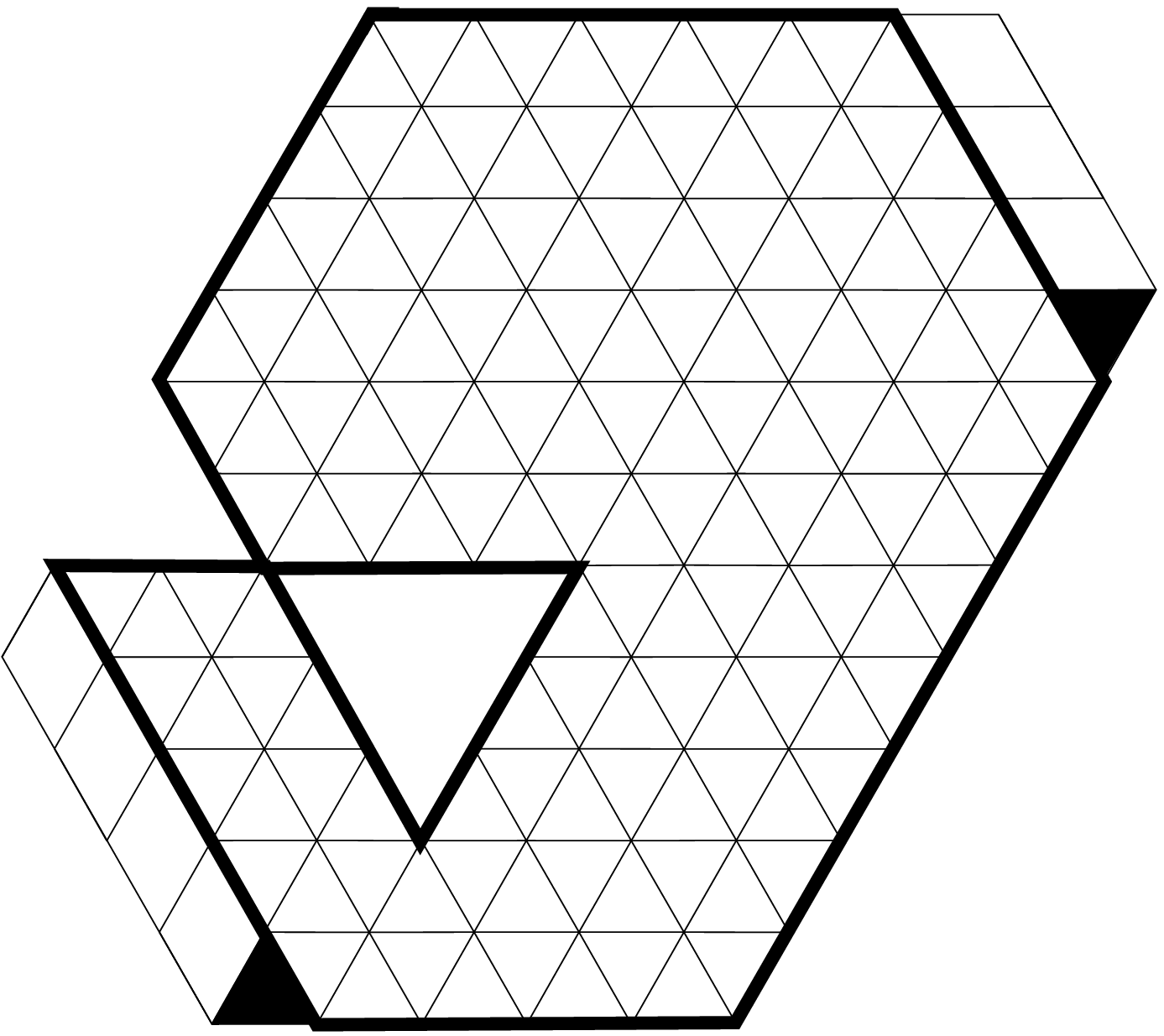}}{\mypic{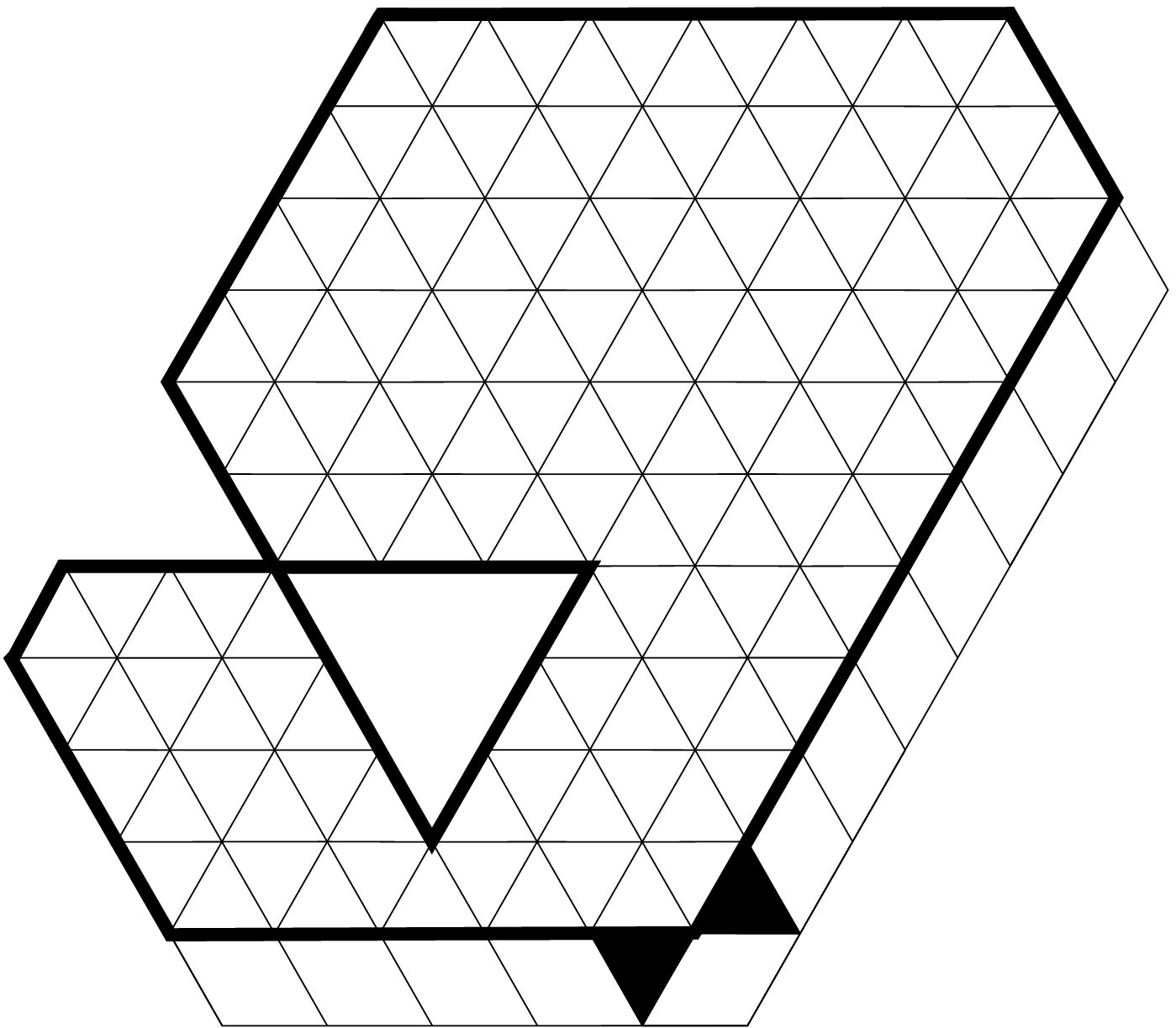}}
\bigskip
\centerline{{\smc Figure {\fcd}.}\ The recurrence for the regions $B_{x,y}(a,b,c,m)$.}

\endinsert

Suppose $x,y,b\geq1$, and assume that (\eca) holds for all values of
the parameters for which the sum of the $x$-parameter, $y$-parameter
and $b$-parameter is less than $x+y+b$. Let $G$ be the planar dual
graph\footnote{By the planar dual graph of a region on the triangular
lattice we understand the graph whose vertices are the unit triangles
inside the region, and whose edges connect vertices corresponding to
unit triangles that share an edge.}
of the region $B_{x,y}(a,b,c,m)$. 
Choose the vertices 
$\alpha$, $\beta$, $\gamma$ and $\delta$
as indicated in Figure~{\fcd},
where $\alpha$ is the leftmost black unit triangle\footnote{By this we understand, of course, the vertex of the planar dual graph corresponding to that unit triangle.}, and $\beta$,
$\gamma$ and $\delta$ are the next black unit triangles as one
moves counterclockwise from $\alpha$ (Figure~{\fcd} corresponds to the
case
$x=3$, $y=1$, $a=4$, $b=1$, $c=3$, $m=2$.)  
Then (\ecb) states that the product of the number of lozenge tilings
of the two regions on top is equal to the product of the number of
lozenge tilings of the two regions in the middle, plus the product of
the number of lozenge tilings of the two regions on the bottom. After
removing the lozenges forced by the unit triangles 
$\alpha$, $\beta$, $\gamma$ and $\delta$ indicated in Figure~{\fcd}, the leftover regions are, in all six instances, magnet bar regions. More precisely, we obtain
$$
\spreadlines{3\jot}
\multline
\kern-10pt
\M(B_{x,y}(a,b,c,m))\M(B_{x-1,y}(a,b-1,c,m))
=
\M(B_{x,y}(a,b-1,c,m))\M(B_{x-1,y}(a,b,c,m))
\\
+
\M(B_{x-1,y+1}(a,b-1,c,m))\M(B_{x,y-1}(a,b,c,m)).
\endmultline
\tag\ecd
$$
Note that all magnet bar regions in the above 
equation except the first one have the sum of their $x$-, $y$- and $b$-parameters strictly less than $x+y+b$. By the induction hypothesis, these five magnet bar regions have their number of lozenge tilings given by (\eca). It is readily checked that substituting these formulas into the above equation one obtains, after simplifications, that $\M(B_{x,y}(a,b,c,m))$ equals precisely the expression on the right-hand side of (\eca). This completes the induction step, and hence the proof. \epf

\mysec{4. Proof of Theorems {\tba} and {\tbb}}

As we have already mentioned, our proof of Theorems~{\tba} and {\tbb}
is also based on Kuo's graphical condensation method. The argument we
use is an induction that proves Theorems~{\tba} and {\tbb}
simultaneously. The previous 
section's result on magnet bar regions is needed for the base cases of our induction.

In our proof it will be essential to be familiar with various distances one can naturally consider within a given $S$-cored hexagon. These are displayed in Figures~{\fdaa} and {\fdab} (note that the unit used is the distance between two consecutive lattice lines).

\topinsert
\centerline{\mypic{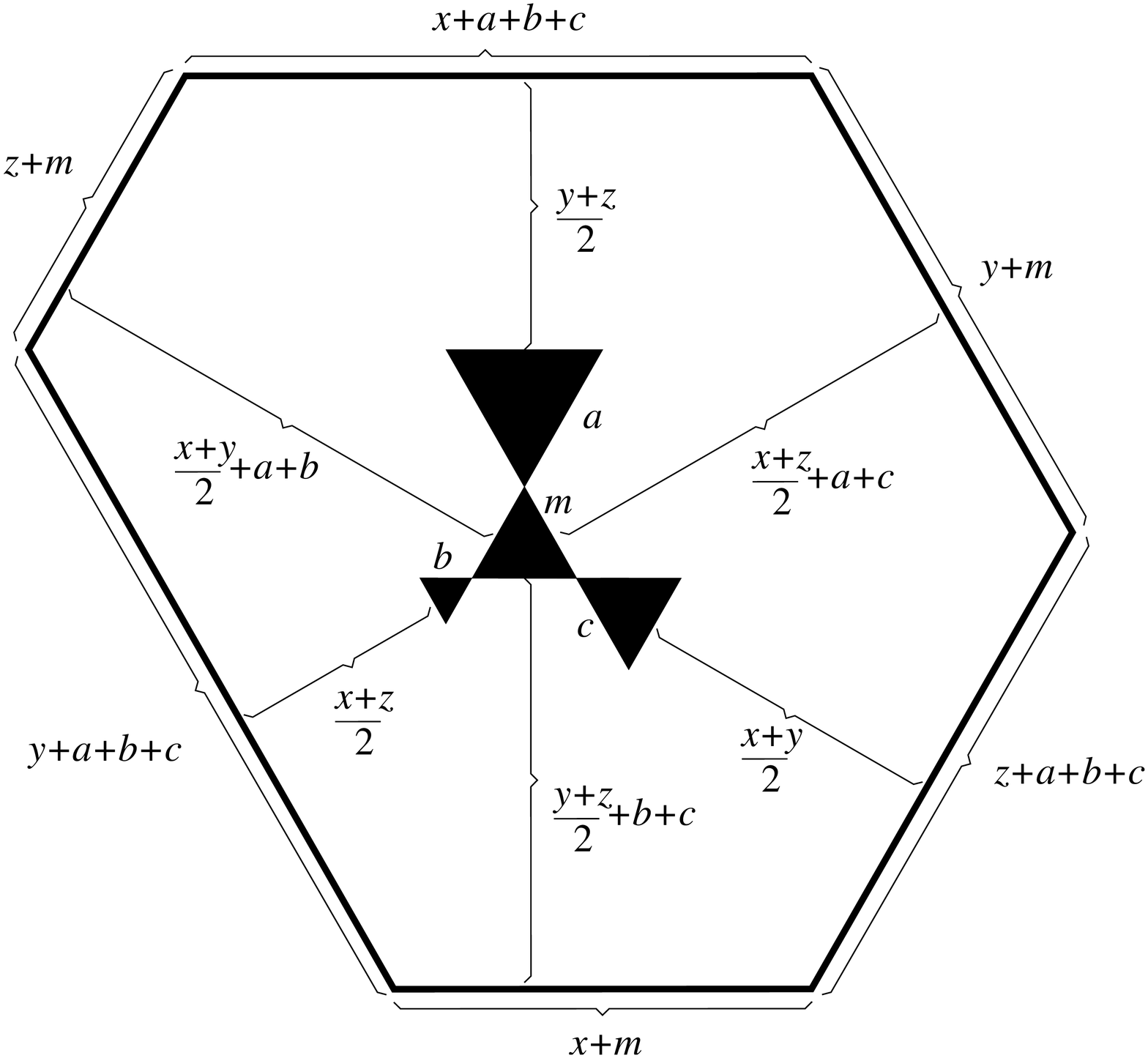}}
\medskip
\centerline{{\smc Figure~{\fdaa}{\rm . Distances within the $S$-cored hexagon $SC_{x,y,z}(a,b,c,m)$}}}
\centerline{\rm when $x$, $y$ and $z$ have the same parity.}
\bigskip
\centerline{\mypic{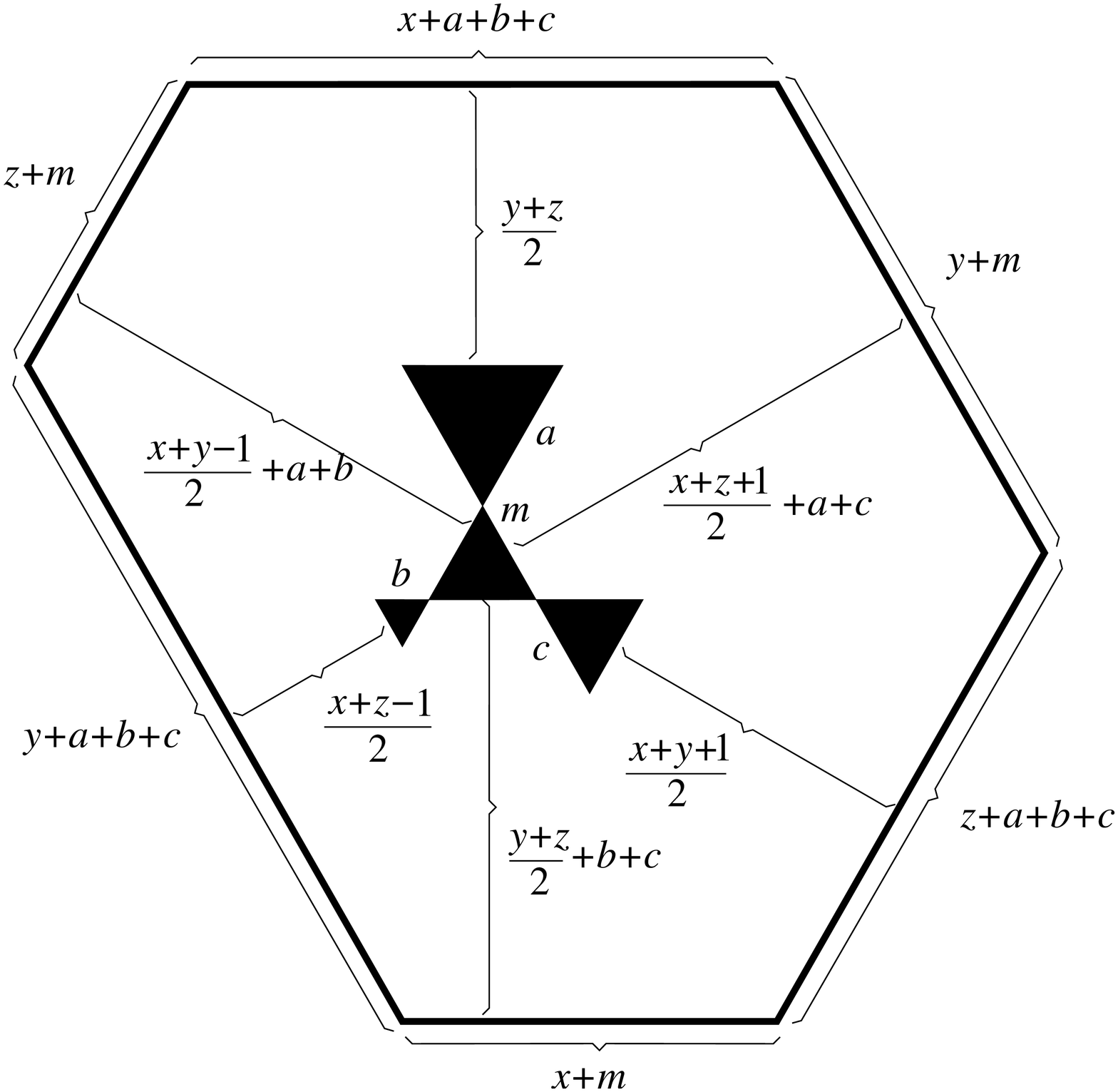}}
\medskip
\centerline{{\smc Figure~{\fdab}{\rm . Distances within an $S$-cored hexagon $SC_{x,y,z}(a,b,c,m)$}}}
\centerline{\rm when $x$ has opposite parity to $y$ and $z$.}
\endinsert

{\it Proof of Theorems {\tba} and {\tbb}.} We prove formulas~(\eba) and (\ebb) by induction on $x+y+z$. The base cases are the instances when $x=0$, $y=0$ or $z=0$.

Consider first the case when $x$, $y$ and $z$ have the same
parity. Due to symmetry, it is enough to verify the instance when
$z=0$. In that case, the $S$-cored hexagon $SC_{x,y,0}(a,b,c,m)$ looks
as illustrated in Figure~{\fda}. Then the hexagon $\H$ cut out from the
$S$-cored hexagon by extending sides of the $a$- and $b$-lobes of the
shamrock core as indicated in Figure~{\fda} (see the dotted lines in
the figure) must be 
{\it internally} tiled in each tiling of $SC_{x,y,0}(a,b,c,m)$. Indeed, this follows by the same argument we used to show that the top hexagonal subregion in Figure~{\fcc} is internally tiled, except we need to apply it now for both dashed-line cuts (the argument applies because the northwestern and southeastern sides of $\H$ have the same length, $m$).

\topinsert
\centerline{\mypic{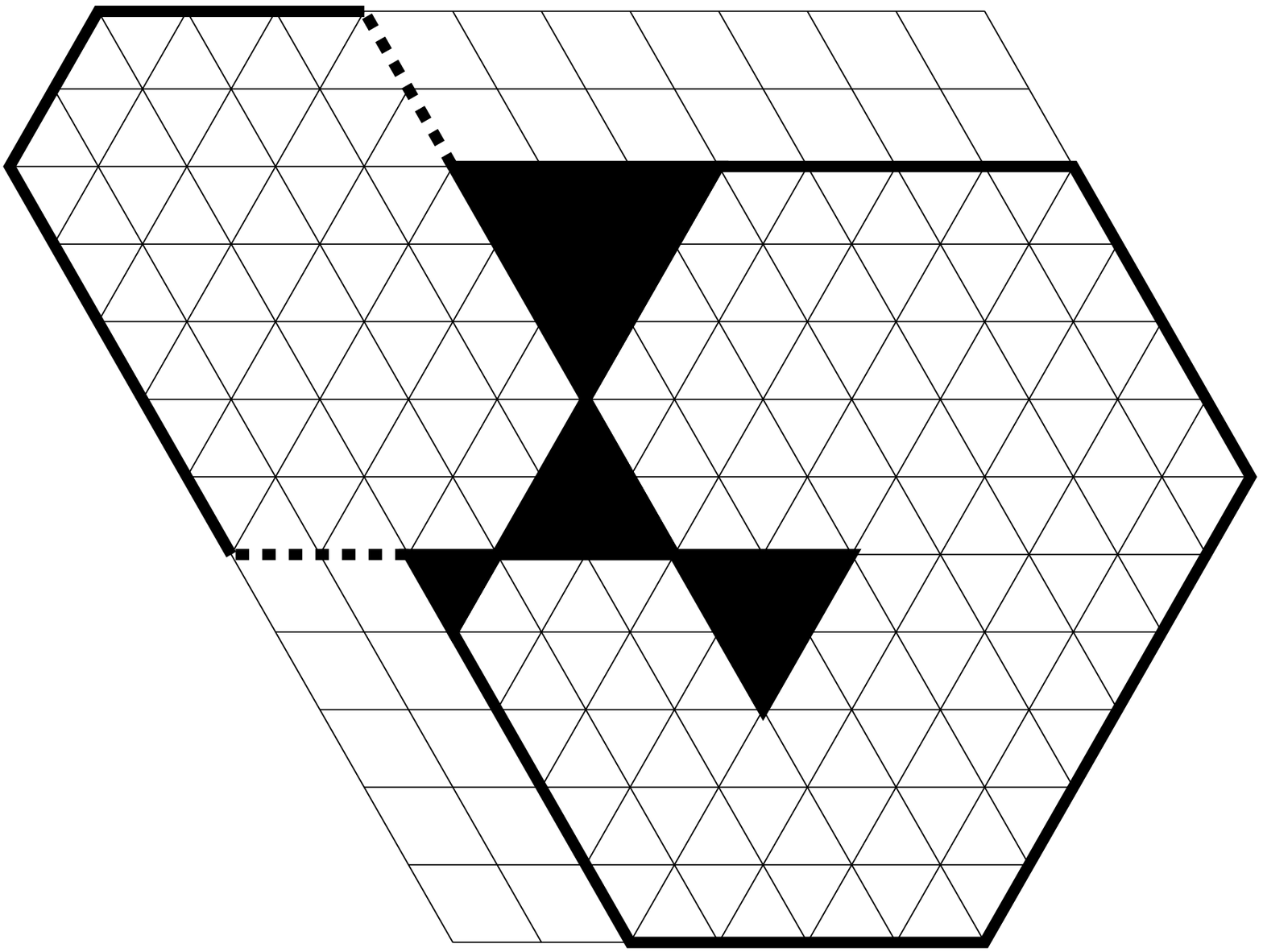}}
\medskip
\centerline{{\smc Figure~{\fda}{\rm . The $S$-cored hexagon $SC_{0,4,4}(3,1,2,2)$.}}}
\endinsert

Since $\H$ is tiled internally, it follows that the rhombus that fits in the top right corner of $SC_{x,y,0}(a,b,c,m)$ and rests on both $\H$ and the $a$-lobe of the shamrock is forced to be tiled as shown. The same is true for the analogous rhombus fitting in the bottom left corner. Note that the region obtained from $SC_{x,y,0}(a,b,c,m)$ by removing $\H$ and these two rhombi is precisely the magnet bar region $B_{\frac{x}{2},\frac{y}{2}}(a,b,c,m)$ (see Figure~{\fdaa}; since we are in the case when $x$, $y$ and $z$ have the same parity and we are assuming $z=0$, both $x$ and $y$ are even). It follows that we have
$$
\spreadlines{3\jot}
\align
\M(SC_{x,y,0}(a,b,c,m))
&=\M(\H)\M(B_{\frac{x}{2},\frac{y}{2}}(a,b,c,m))
\\
&=P\left(m,\frac{x}{2}+b,\frac{y}{2}+a\right)
\M(B_{\frac{x}{2},\frac{y}{2}}(a,b,c,m))
\tag\eda
\endalign
$$
(for the second equality we used Figure~{\fdaa} to read off the
side-lengths of $\H$). 
Substituting the values of $P(m,\frac{x}{2}+b,\frac{y}{2}+a)$ and $\M(B_{\frac{x}{2},\frac{y}{2}}(a,b,c,m))$ given by formulas~(\eaa) and (\eca) in the right-hand side above, it is apparent that it becomes precisely the $z=0$ specialization of the expression on the right-hand side of (\eba). This checks the base case of the induction when $x$, $y$ and $z$ have the same parity.

When $x$, $y$ and $z$ have mixed parities, we can assume without loss
of generality that $x$ has parity 
opposite to the parities of $y$ and $z$.
There are now two inequivalent base cases to check, namely $x=0$ and $z=0$ ($y=0$ is equivalent to the latter). Consider first the case when $x=0$ (note that our assumptions imply then that $y$ and $z$ are odd). Then the $S$-cored hexagon $SC_{0,y,z}(a,b,c,m)$ looks as illustrated in Figure~{\fdb}.

\topinsert
\twoline{\mypic{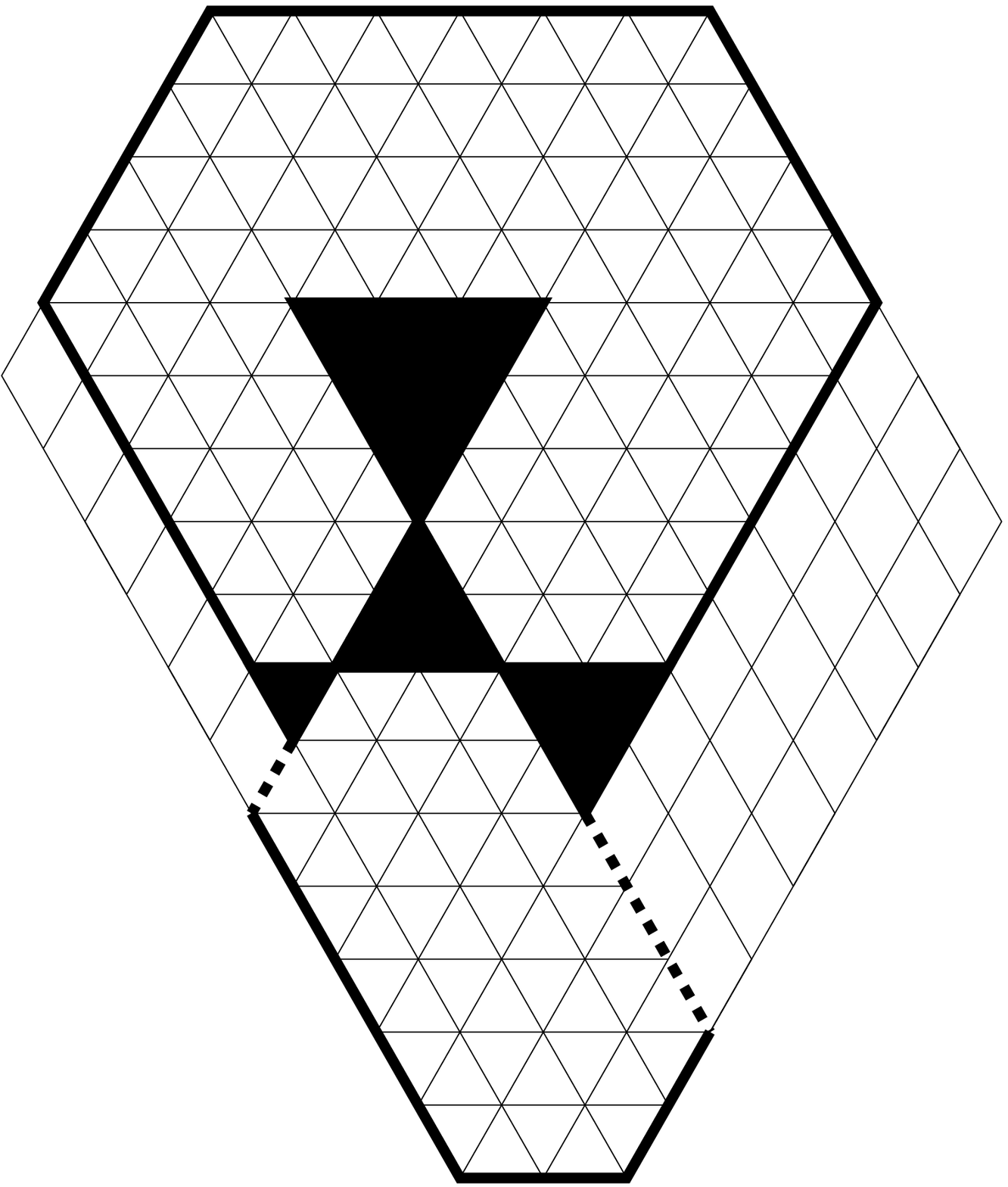}}{\mypic{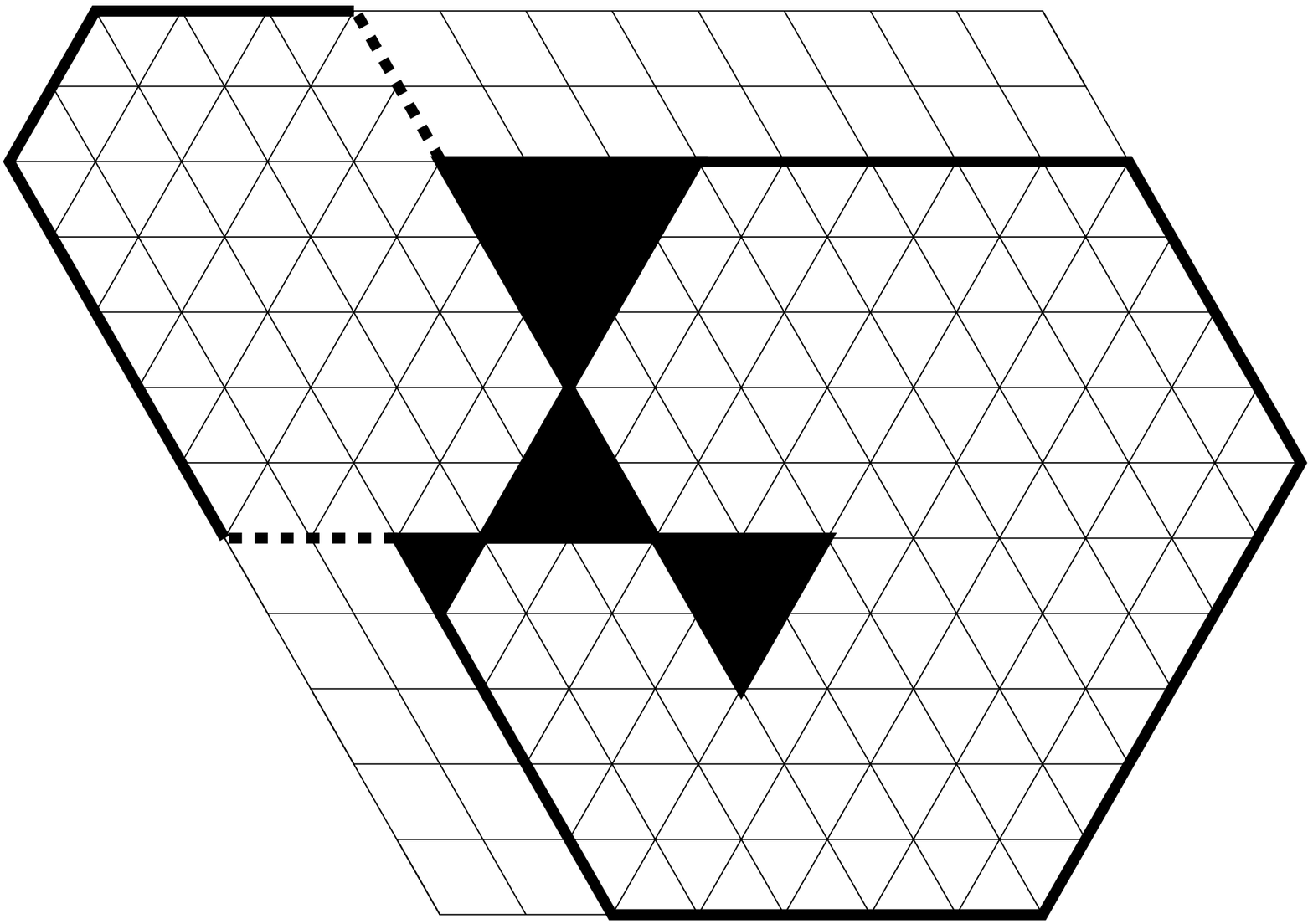}}
\medskip
\twoline{{\smc Figure {\fdb}.}\ The region $SC_{0,5,3}(3,1,2,2)$.}
{{\smc Figure {\fdc}.}\ The region $SC_{5,4,0}(3,1,2,2)$.}
\endinsert

Note that, by the argument we used in the paragraph before (\eda), the hexagonal subregion determined by the dashed lines in Figure~{\fdb} must be internally matched in any lozenge tiling of $SC_{0,y,z}(a,b,c,m)$. This implies that the two rhombic regions at the left and right corners of the $S$-cored hexagon whose lozenge tilings are indicated in Figure~{\fdb} are tiled as shown in each tiling of $SC_{0,y,z}(a,b,c,m)$. The region left over from the $S$-cored hexagon after removing these two rhombic regions is the disjoint union of the hexagonal subregion mentioned earlier in this paragraph, and a magnet bar region. It follows from the distances indicated in Figure~{\fdab} that the three values for the hexagon side-lengths are $m$, $\frac{y+1}{2}+c$ and $\frac{z-1}{2}+b$, while the magnet bar region is $B_{\frac{y-1}{2},\frac{z+1}{2}}(b,c,a,m)$. Thus we have
$$
\M(SC_{0,y,z}(a,b,c,m))=
P\left(m,\frac{y+1}{2}+c,\frac{z-1}{2}+b\right)
\M(B_{\frac{y-1}{2},\frac{z+1}{2}}(b,c,a,m)).
\tag\edb
$$
Using formulas~(\eaa) and (\eca) one readily sees that the above equation agrees with the $x=0$ specialization of (\ebb).

The remaining base case to check is $z=0$, $x$ odd and $y$ even. This follows by the same arguments as the other two base cases. Figure~{\fdc} shows how $\M(SC_{x,y,0}(a,b,c,m))$ is the product of the number of tilings of a hexagon and a magnet bar region. Their precise dimensions can be deduced from Figure~{\fdab}. One obtains
$$
\M(SC_{x,y,0}(a,b,c,m))=
P\left(m,\frac{x-1}{2}+b,\frac{y}{2}+a\right)
\M(B_{\frac{x+1}{2},\frac{y}{2}}(a,b,c,m)).
\tag\edc
$$
Using (\eaa) and (\eca), this is again readily seen to agree with the $z=0$ specialization of (\ebb). This concludes the verification of all the base cases we need for our induction.

\topinsert

\twoline{\mypic{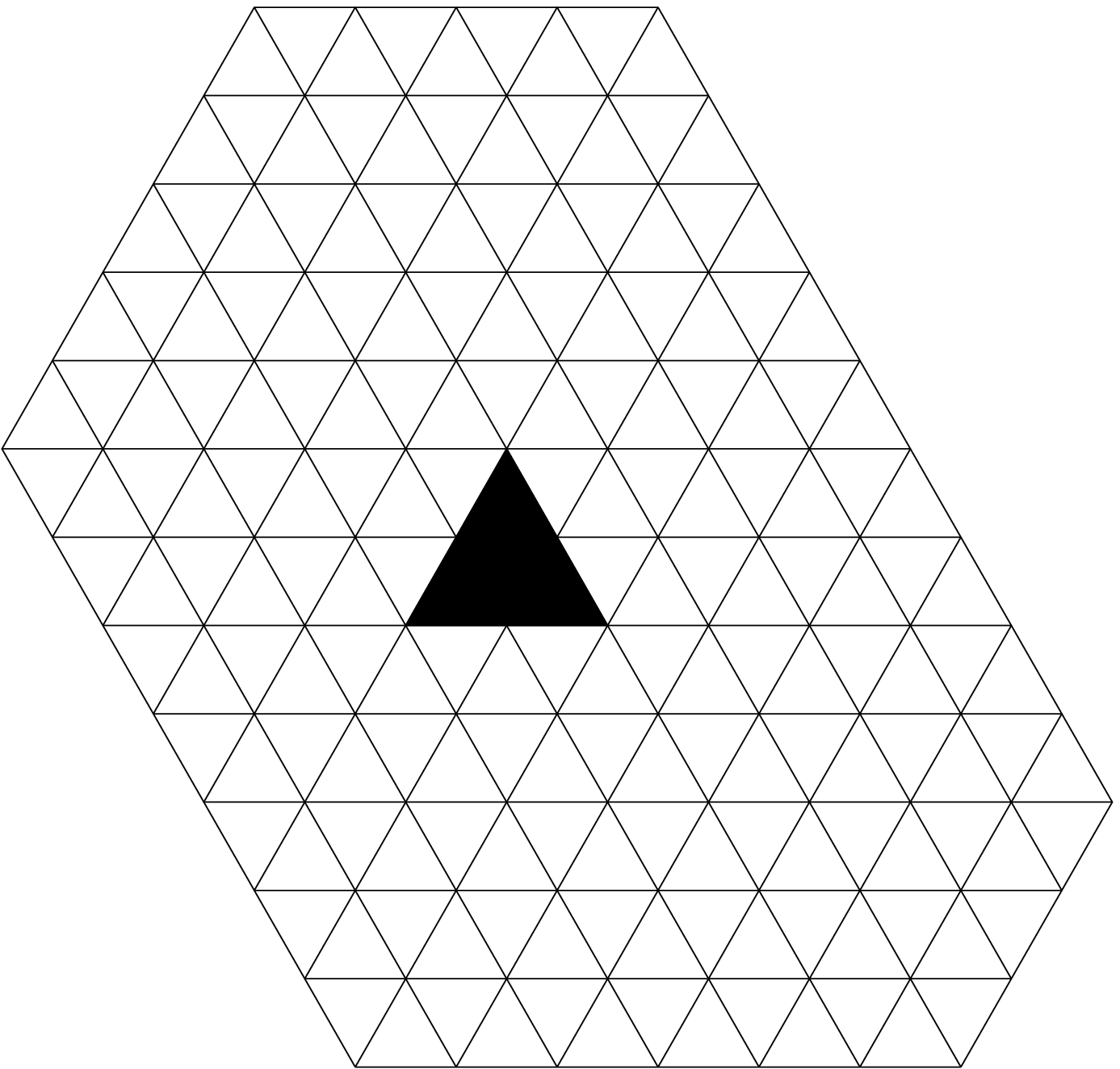}}{\mypic{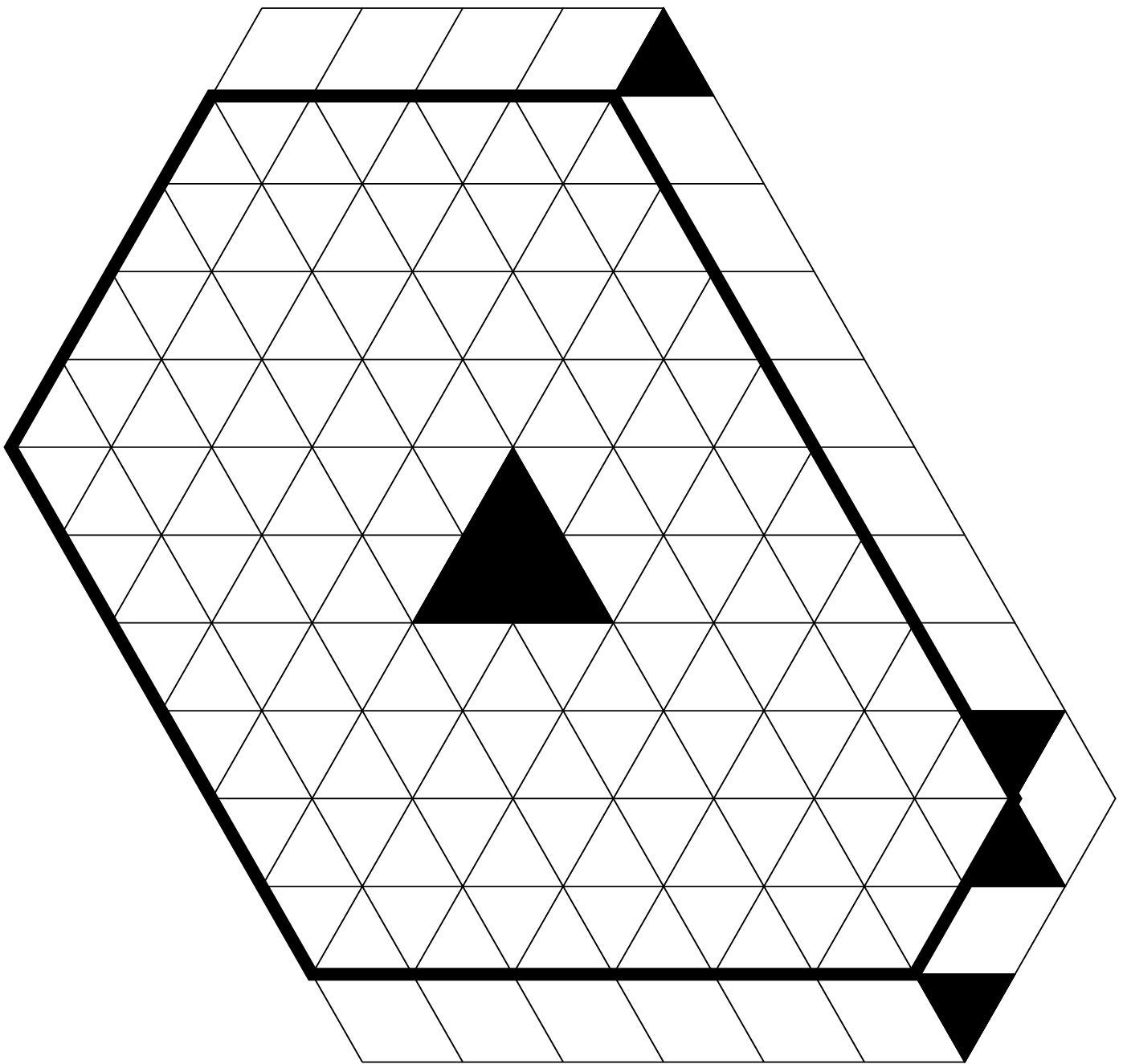}}
\bigskip

\twoline{\mypic{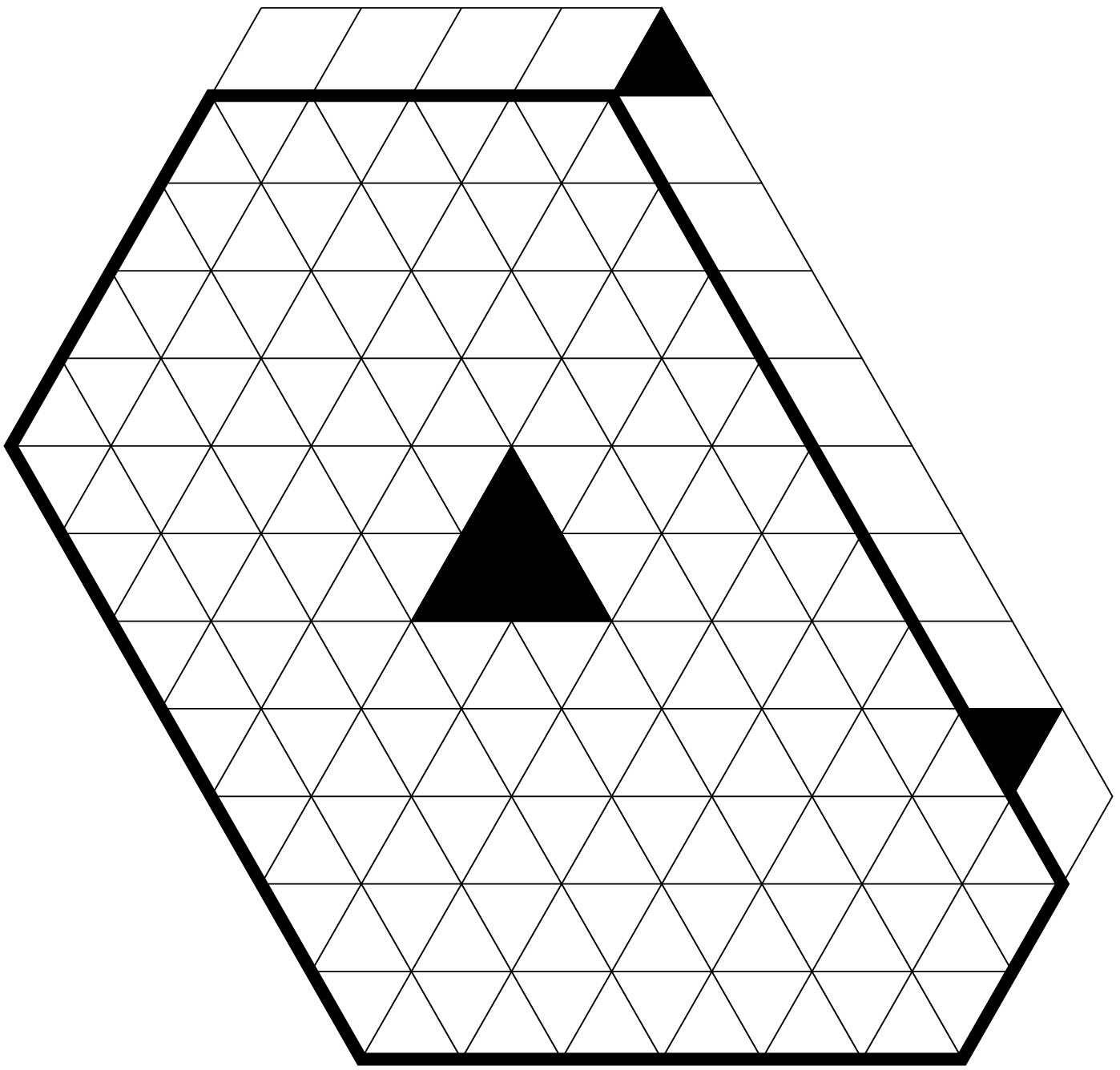}}{\mypic{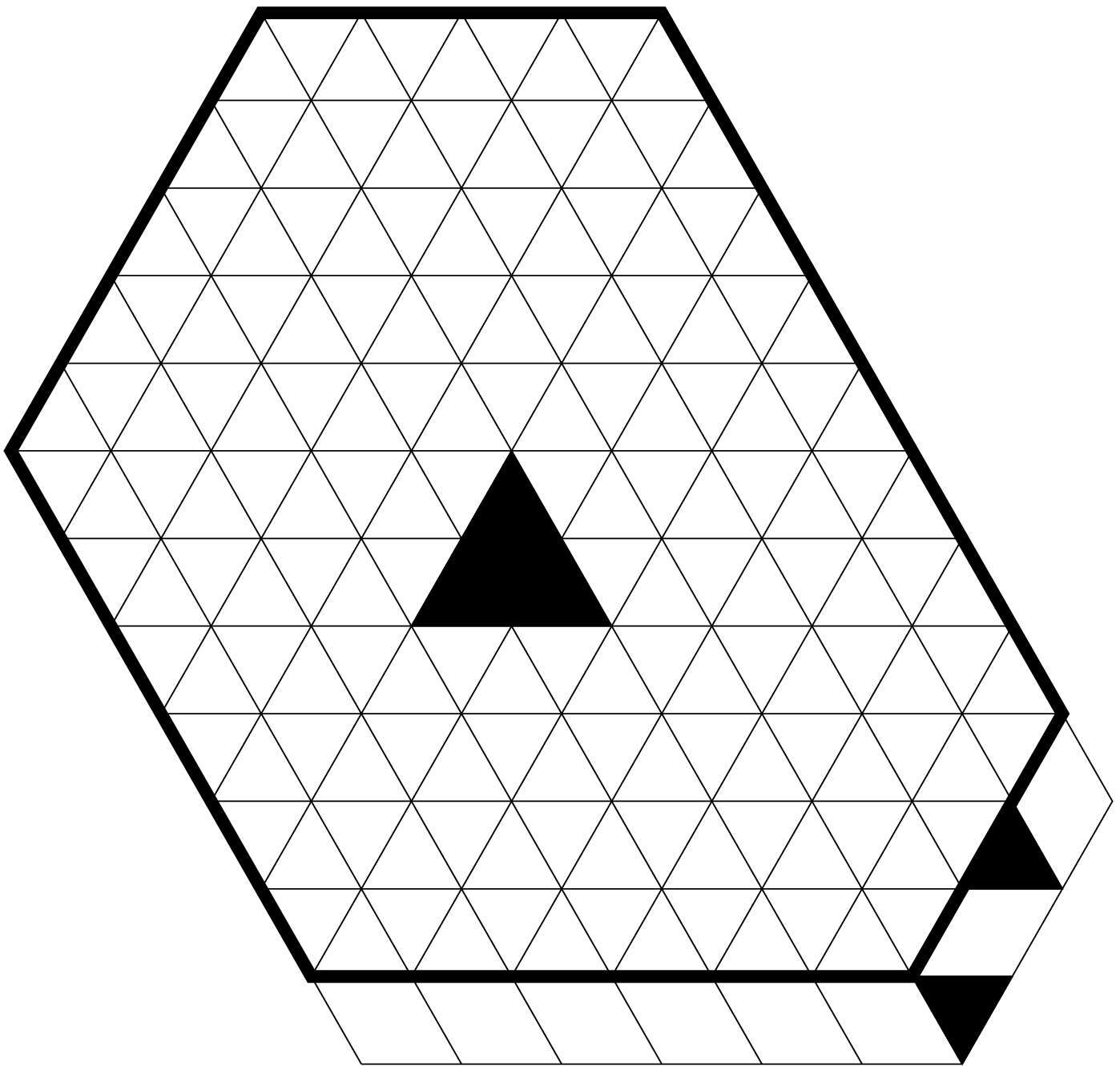}}
\bigskip

\twoline{\mypic{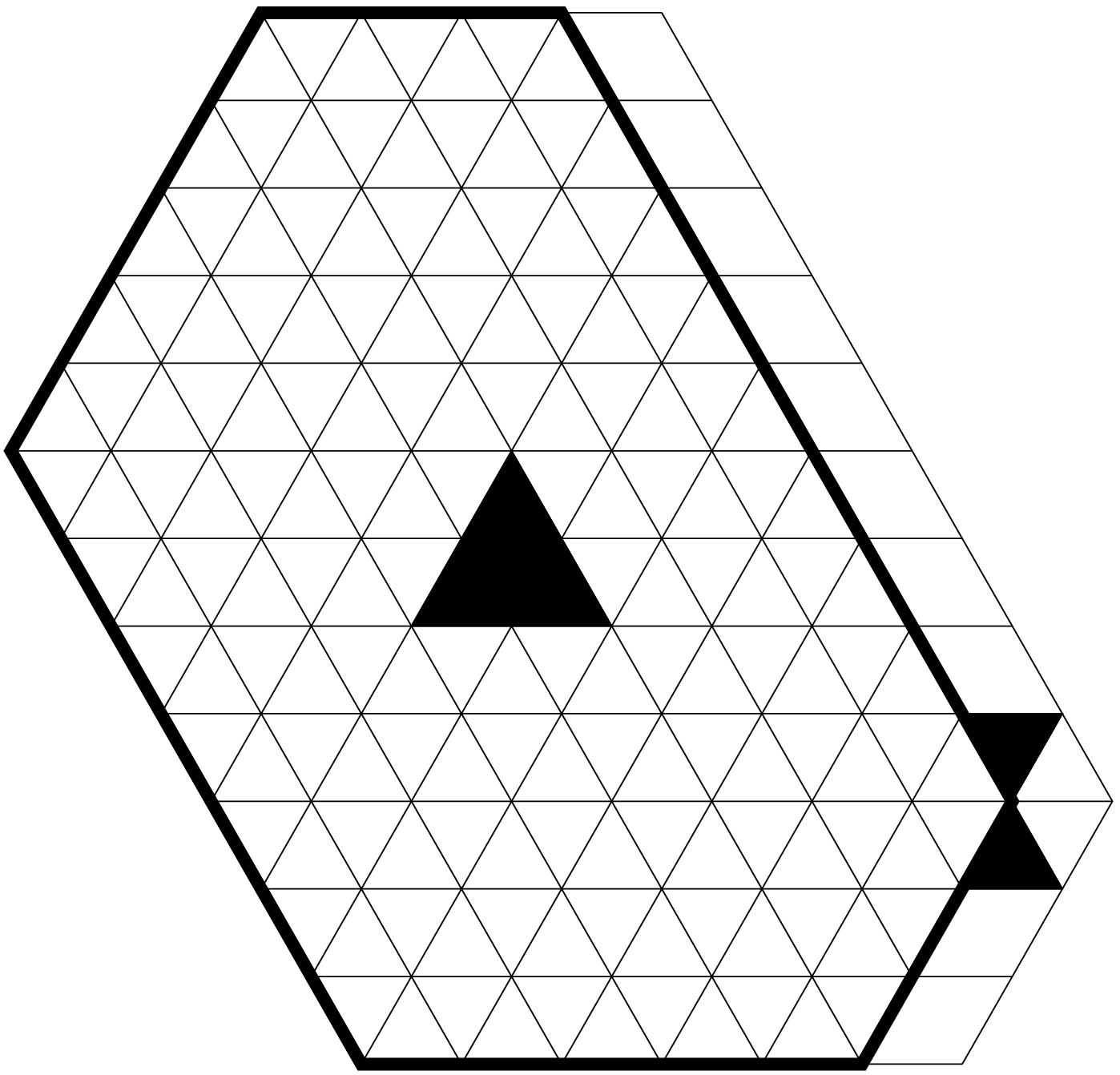}}{\mypic{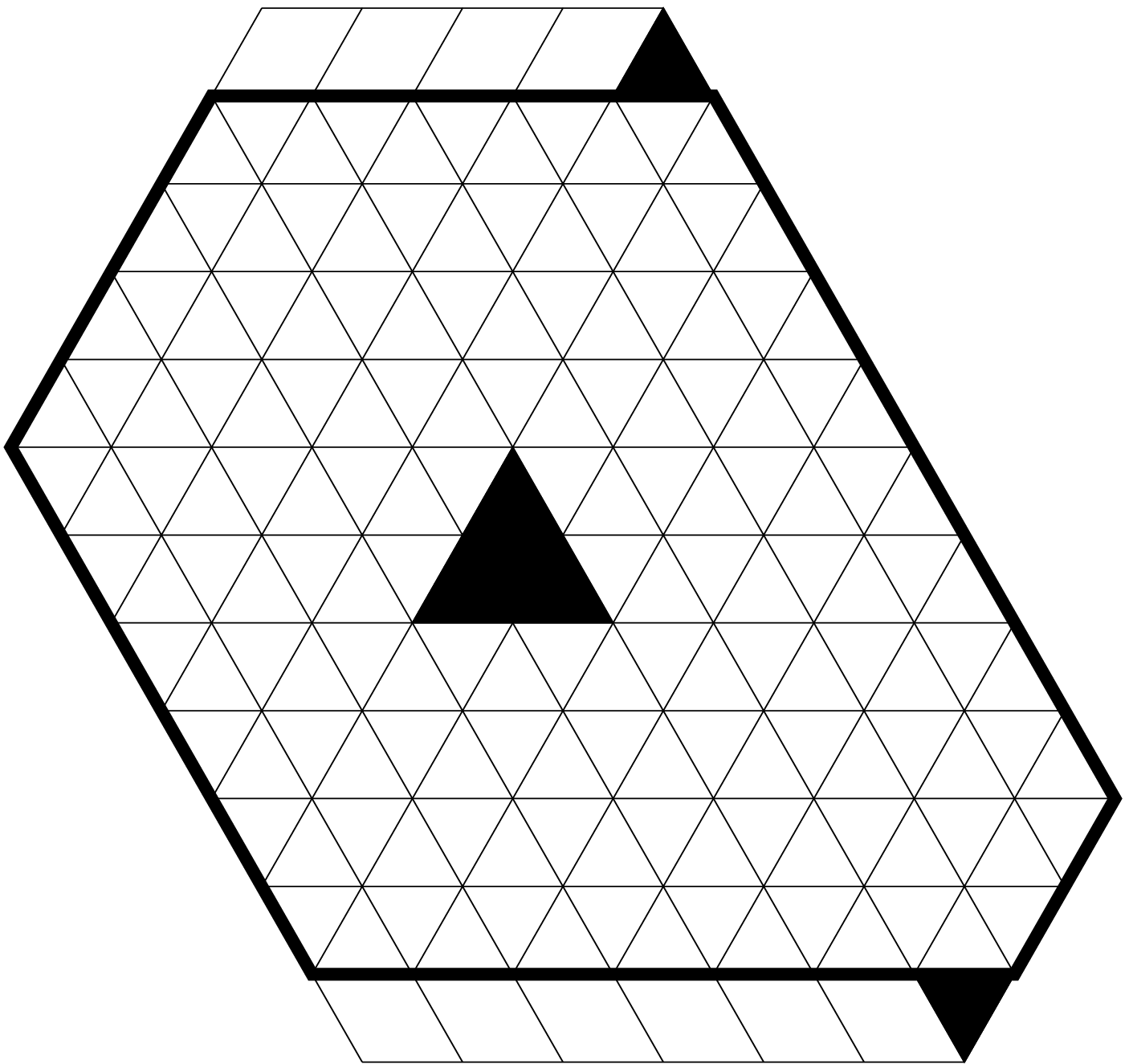}}
\bigskip
\centerline{{\smc Figure {\fdd}.}\ The recurrence for the regions $SC_{x,y,z}(0,0,0,m)$, $x$ of opposite parity to $y$ and $z$.}

\endinsert

The induction step is based on Kuo's graphical condensation (stated in Theorem~{\tcb}). More precisely, we use graphical condensation to obtain recurrences for the number of lozenge tilings of the $S$-cored hexagons $SC_{x,y,z}(a,b,c,m)$, valid for $x,y,z\geq1$, and then verify that the right-hand side of the claimed formulas~(\eba) and (\ebb) satisfy the same recurrences.

We will need two different recurrences, one for the case when $x$, $y$
and $z$ have the same parity, and one for the mixed parity case. Both
will follow by applying graphical condensation in the same fashion as
in the proof of Theorem~{\tca}, 
namely by choosing the four removed unit triangles to be along two sides of the outer boundary of the $S$-cored hexagons, in the pattern shown in Figure~{\fcd}. In particular, only the outer boundary of the $S$-cored hexagons will change in the resulting subregions, while the shamrock core remains intact. Due to this, it suffices to discuss how we obtain our recurrences in the case when $a=b=c=0$.

Consider first the case when $x$, $y$ and $z$ have mixed parities (an instance of this is illustrated in Figure~{\fdd}). Without loss of generality we may assume that the parity of $x$ is different from that of $y$ and $z$. When applying Kuo's graphical condensation as described in the previous paragraph, it is of crucial importance in which corner of the outer hexagon one places the pattern of the four removed unit triangles. Since in $SC_{x,y,z}(0,0,0,m)$ (which is the same as the cored hexagon $C_{x,y,z}(m)$ of \cite{\cekz}) the core is just to the left of the true central position, the correct choice will be to place the pattern of the four removed unit triangles in the right corner of the outer hexagon (see Figure~{\fdd}).

Let $G$ be the planar dual graph of the region
$SC_{x,y,z}(0,0,0,m)$. Choose the vertices 
$\alpha$, $\beta$, $\gamma$ and $\delta$ as
indicated in Figure~{\fdd}, where $\alpha$ is
the bottom black unit
triangle, and $\beta$, $\gamma$ and $\delta$ are
the next black unit triangles as
one moves upwards (Figure~{\fdd} corresponds to the case 
$x=4$, $y=7$, $z=3$). Then (\ecb) states that the product of the
number of lozenge tilings of the two regions on top is equal to the
product of the number of lozenge tilings of the two regions in the
middle, plus the product of the number of lozenge tilings of the two
regions on
the bottom. After removing the lozenges forced by the black
unit triangles indicated in Figure~{\fdd}, the leftover regions are,
in all six instances, $S$-cored hexagons. (Note that this would not be
the case if we applied graphical condensation with the pattern of
removed triangles placed in a different corner --- some of the
resulting positions of the core would not be central!) The precise
parameters of these $S$-cored hexagons can be visually extracted from
Figure~{\fdd} (this is easier to do in the case 
$a=b=c=0$; this is why we have reduced the proof of the recurrences to this case). 

Indeed, the leftover region on the top right of Figure~{\fdd} is $SC_{x,y-1,z-1}(0,0,0,m)$. The central left region, after clockwise rotation by $120^\circ$ followed by a reflection across the vertical, becomes $SC_{y,x,z-1}(0,0,0,m)$. The central right region, after counterclockwise rotation by $120^\circ$ followed by a reflection across the vertical, becomes $SC_{z,y-1,x}(0,0,0,m)$. Similarly, the bottom two leftover regions are $SC_{x-1,y,z}(0,0,0,m)$ and $SC_{x+1,y-1,z-1}(0,0,0,m)$, respectively. We obtain 
$$
\spreadlines{3\jot}
\multline
\M(SC_{x,y,z}(0,0,0,m))\M(SC_{x,y-1,z-1}(0,0,0,m))
\\
=
\M(SC_{y,x,z-1}(0,0,0,m))\M(SC_{z,y-1,x}(0,0,0,m))
\\
+
\M(SC_{x-1,y,z}(0,0,0,m))\M(SC_{x+1,y-1,z-1}(0,0,0,m)).
\endmultline
\tag\edd
$$
As explained above, this yields a recurrence for any nonnegative
values for the sizes $a$, $b$ and $c$ of the lobes of the shamrock
core. Taking into account the two rotations followed by reflections
that we needed to consider when converting Figure~{\fdd} into
recurrence~(\edd), one sees that the resulting recurrence for
arbitrary 
$a$, $b$, $c$ is
$$
\spreadlines{3\jot}
\multline
\M(SC_{x,y,z}(a,b,c,m))\M(SC_{x,y-1,z-1}(a,b,c,m))
\\
=
\M(SC_{y,x,z-1}(b,a,c,m))\M(SC_{z,y-1,x}(c,b,a,m))
\kern4cm
\\
+
\M(SC_{x-1,y,z}(a,b,c,m))\M(SC_{x+1,y-1,z-1}(a,b,c,m)),
\\
x,y,z\geq1,\ \text{\rm $x$ of opposite parity to $y$ and $z$}.
\endmultline
\tag\ede
$$

\topinsert

\twoline{\mypic{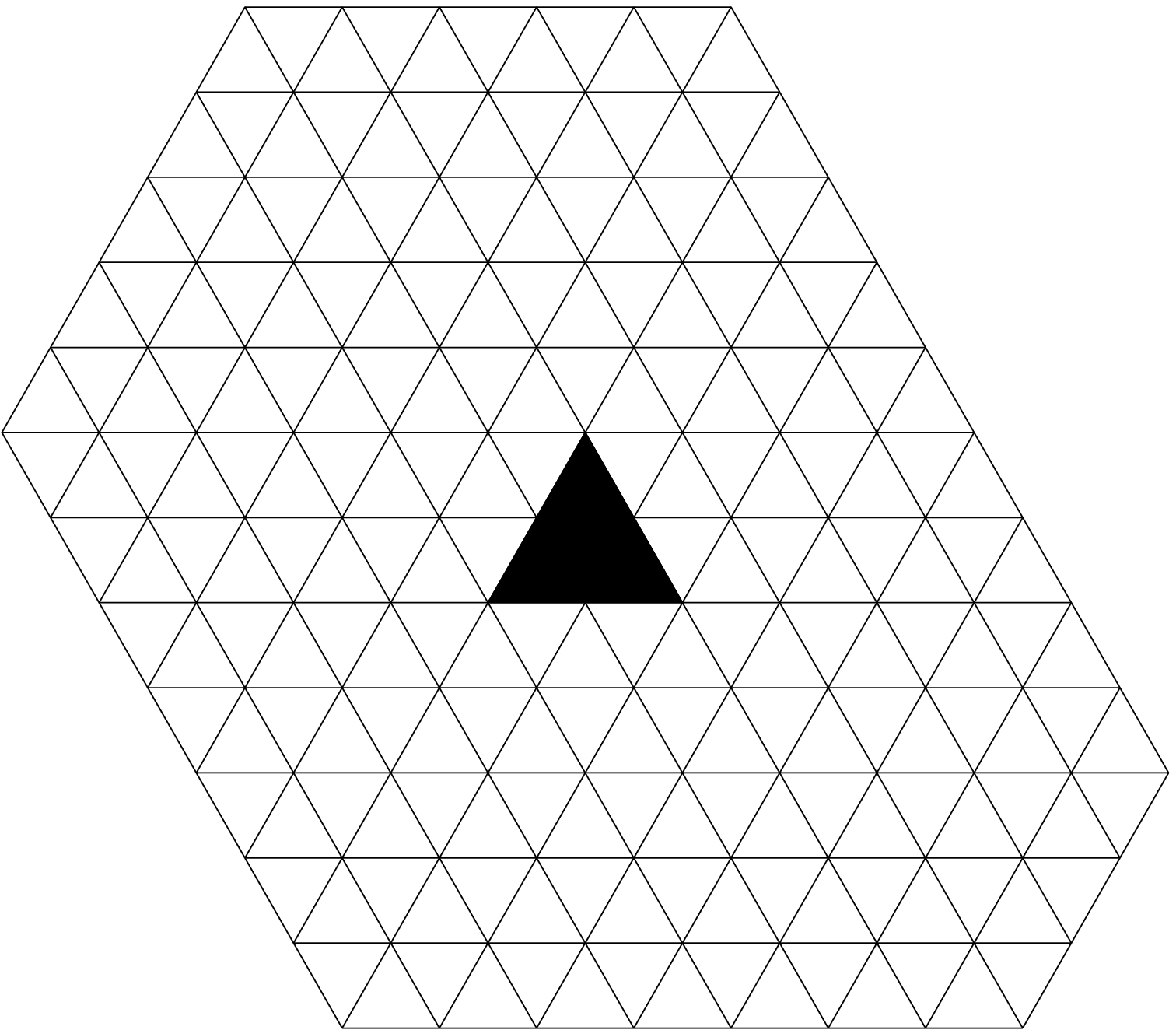}}{\mypic{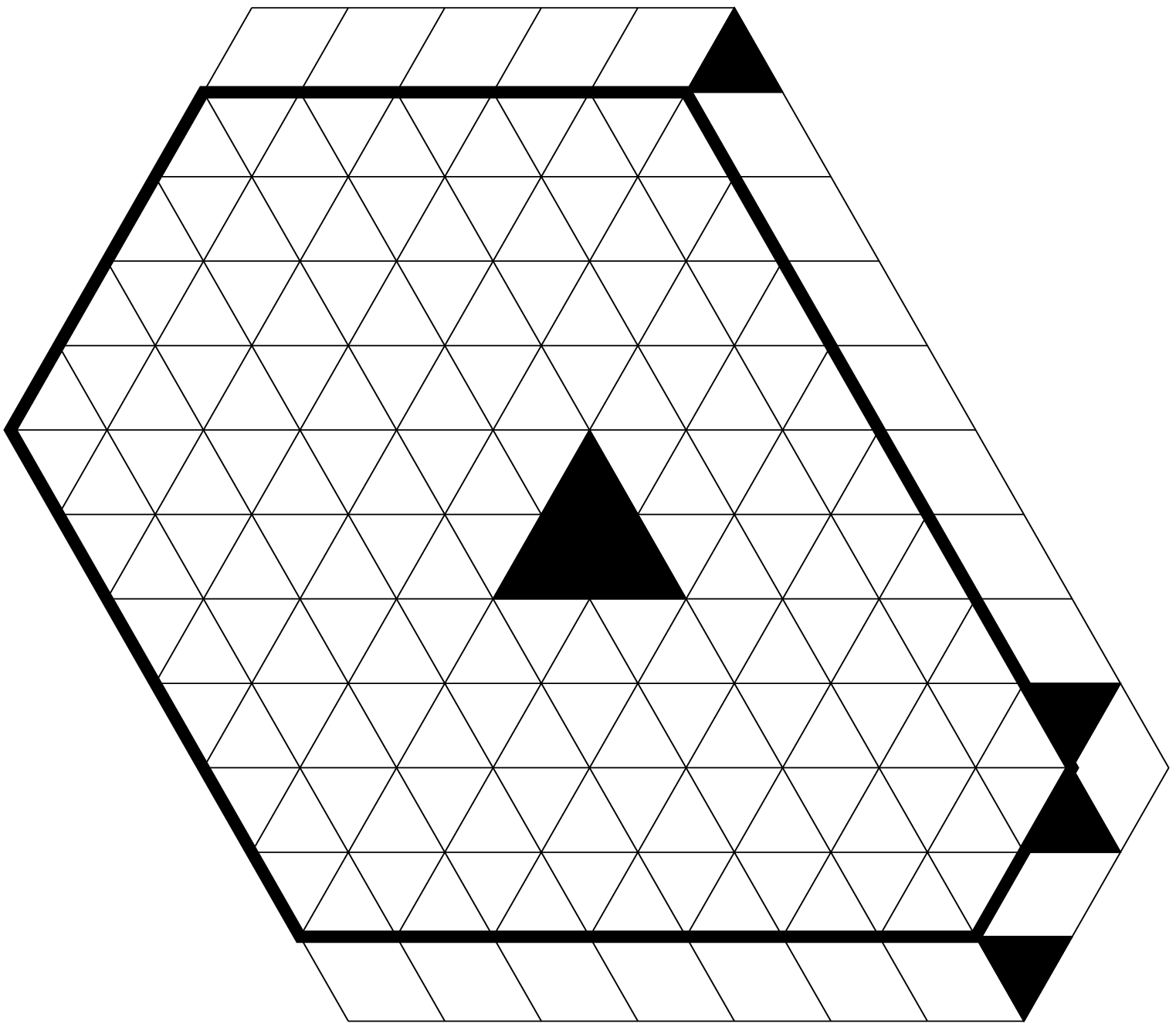}}
\bigskip

\twoline{\mypic{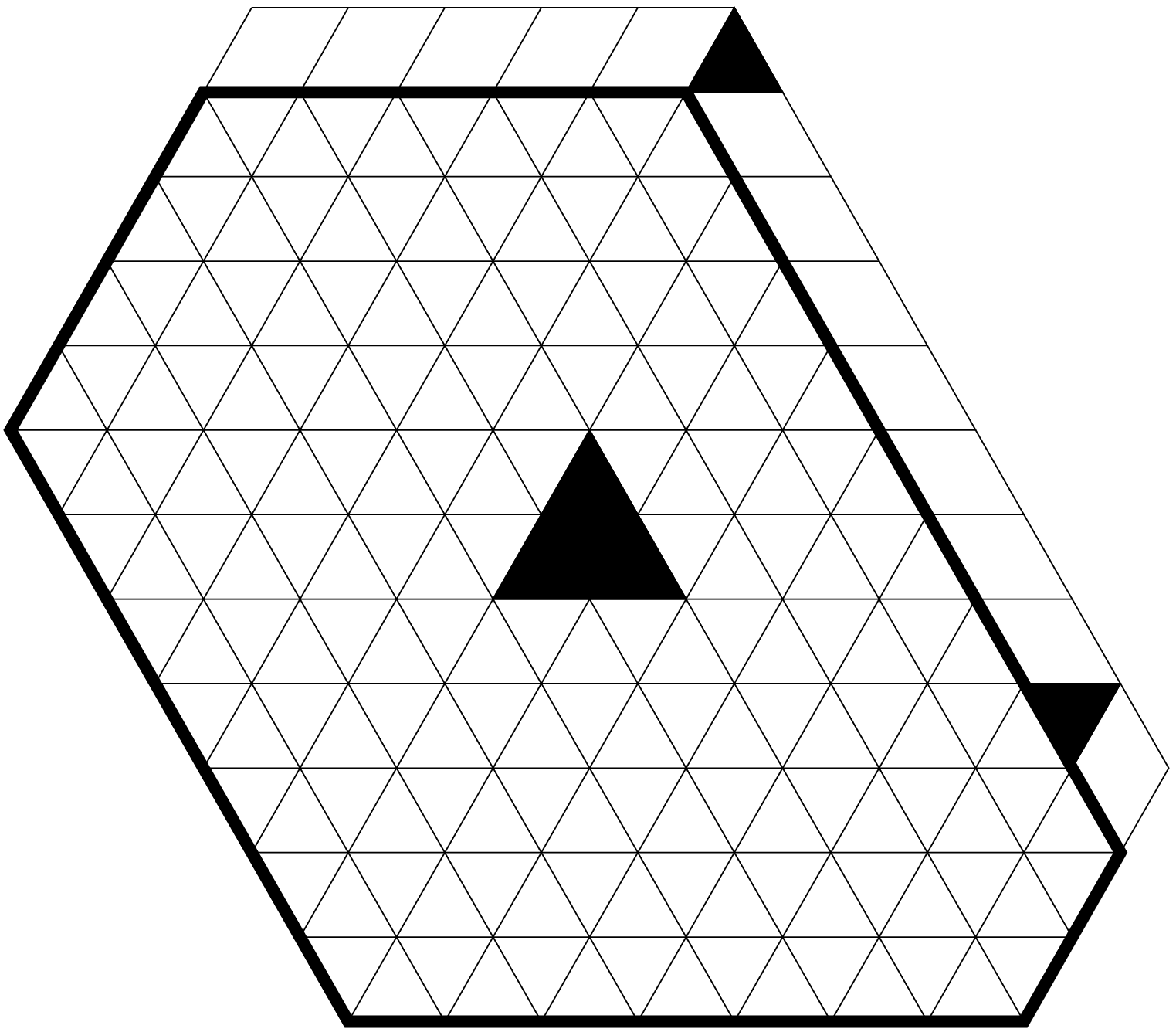}}{\mypic{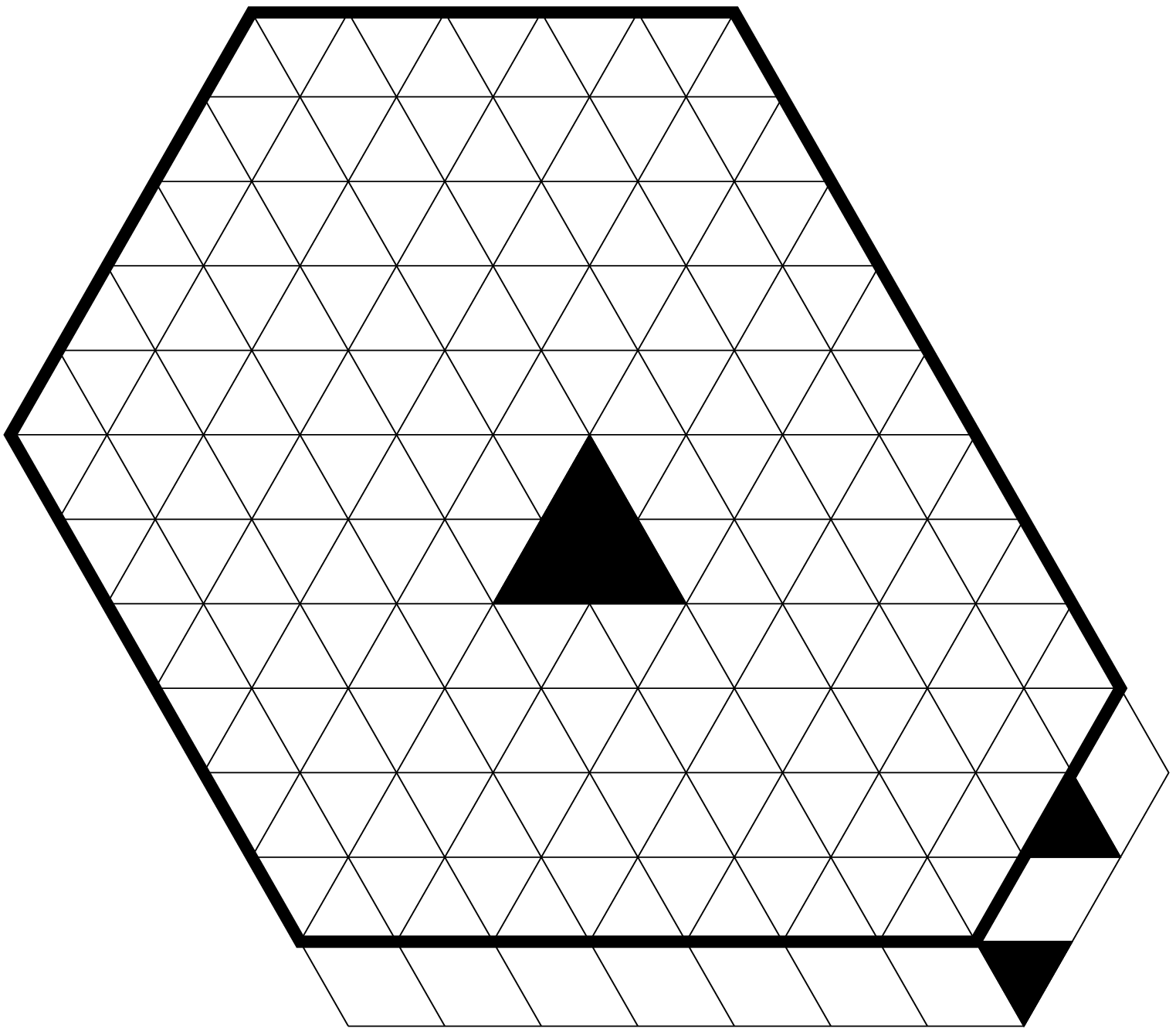}}
\bigskip

\twoline{\mypic{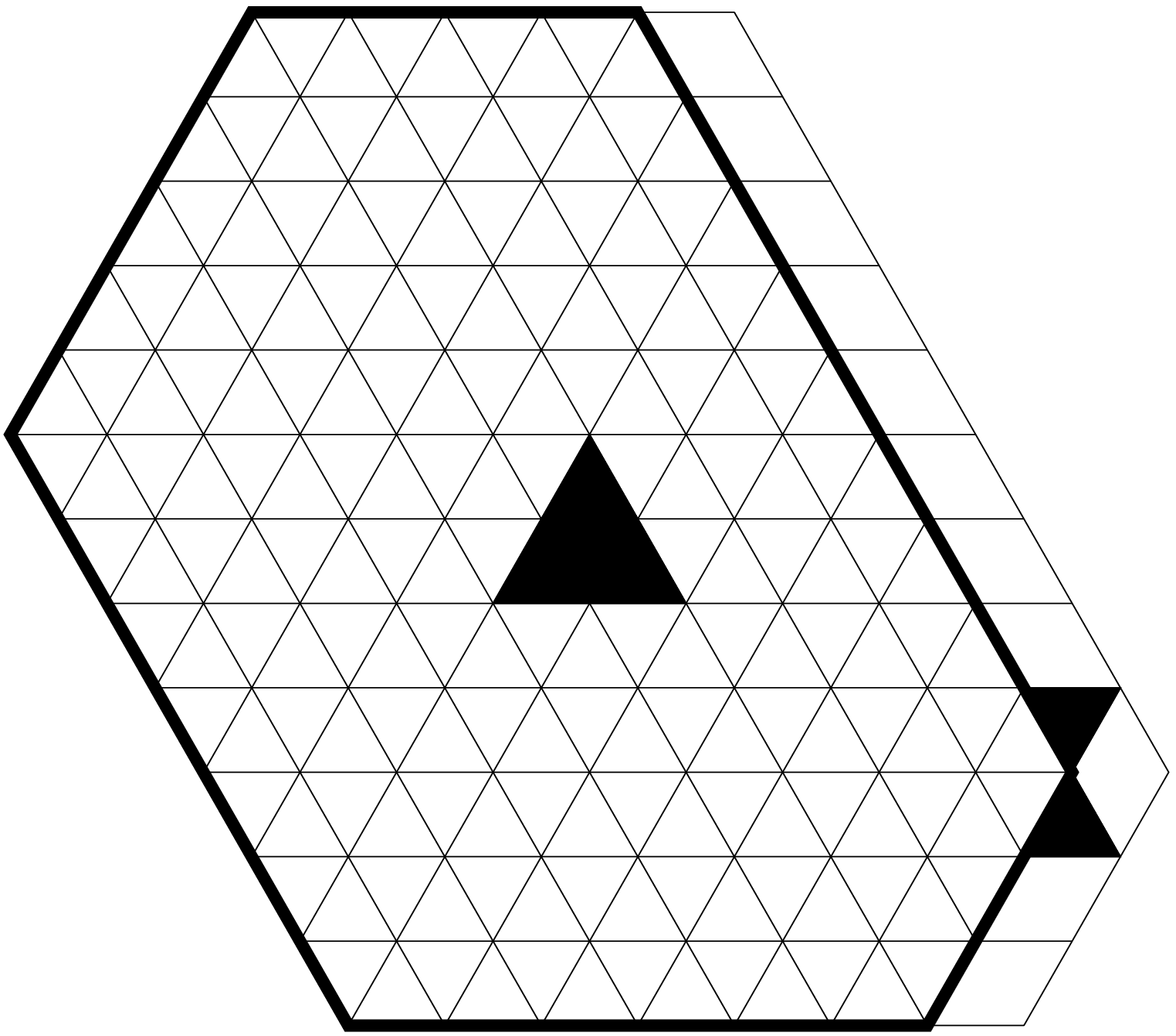}}{\mypic{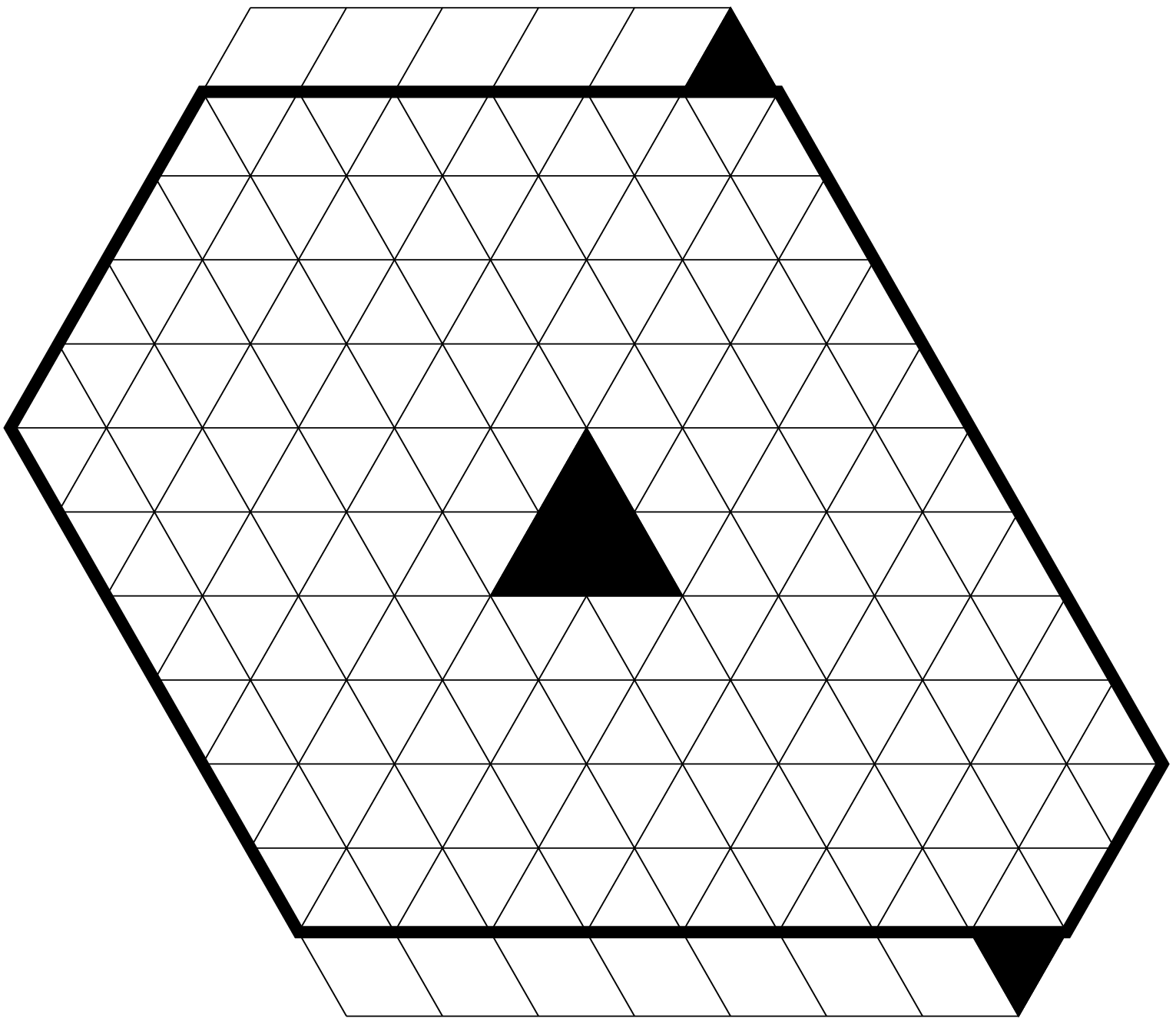}}
\bigskip
\centerline{{\smc Figure {\fde}.}\ The recurrence for the regions $SC_{x,y,z}(0,0,0,m)$, $x$, $y$ and $z$ of same parity.}

\endinsert

In the remaining case when $x$, $y$ and $z$ have the same parity,
there is no ``singled out'' corner of the outer hexagon in which to
place the pattern of removed unit triangles when applying condensation
--- any of the six choices gives a recurrence for $S$-cored
hexagons. In order to be closer to the previous case, we choose again
to fit this pattern in the rightmost corner; see Figure~{\fde}, which
illustrates the case 
$x=5$, $y=7$, $z=3$. Then (\ecb) states that the product of the
number of lozenge tilings of the two regions on top is equal to the
product of the number of lozenge tilings of the two regions in the
middle, plus the product of the number of lozenge tilings of the two
regions on
the bottom. After removing the lozenges forced by the black unit triangles indicated in Figure~{\fde}, all the leftover regions are $S$-cored hexagons. Their precise parameters can be easily deduced from Figure~{\fde}. 
Indeed, the region on the top right, after reflection in the vertical, becomes $SC_{x,z-1,y-1}(0,0,0,m)$. After counterclockwise rotation by $120^\circ$, the left central region becomes $SC_{z-1,y,x}(0,0,0,m)$. The right central region is transformed into $SC_{y-1,z,x}(0,0,0,m)$ by clockwise rotation by $120^\circ$. The bottom left region, when reflected across the vertical, becomes $SC_{x-1,z,y}(0,0,0,m)$. Finally, the bottom right region is precisely $SC_{x+1,y-1,z-1}(0,0,0,m)$. We obtain 
$$
\spreadlines{3\jot}
\multline
\M(SC_{x,y,z}(0,0,0,m))\M(SC_{x,z-1,y-1}(0,0,0,m))
\\
=
\M(SC_{z-1,y,x}(0,0,0,m))\M(SC_{y-1,z,x}(0,0,0,m))
\\
+
\M(SC_{x-1,z,y}(0,0,0,m))\M(SC_{x+1,y-1,z-1}(0,0,0,m)).
\endmultline
\tag\edf
$$
By the above discussion, taking into account the symmetries we needed to apply to the regions in Figure~{\fde} in the previous paragraph, this extends to arbitrary values of $a$, $b$ and $c$. The resulting recurrence is
$$
\spreadlines{3\jot}
\multline
\M(SC_{x,y,z}(a,b,c,m))\M(SC_{x,z-1,y-1}(a,c,b,m))
\\
=
\M(SC_{z-1,x,y}(c,a,b,m))\M(SC_{y-1,z,x}(b,c,a,m))
\kern4cm
\\
+
\M(SC_{x-1,z,y}(a,c,b,m))\M(SC_{x+1,y-1,z-1}(a,b,c,m)),
\\
x,y,z\geq1,\ \text{\rm $x$, $y$ and $z$ of same parity}.
\endmultline
\tag\edg
$$
Let $x,y,z\geq1$ and $a$, $b$, $c$ and $m$ be fixed nonnegative integers, and assume that (\eba) and (\ebb) hold for all values of the parameters for which the sum of the $x$-parameter, $y$-parameter and $z$-parameter is less than $x+y+z$. We need to show that (\eba) and (\ebb) hold also for $x$, $y$, $z$, $a$, $b$, $c$ and $m$.

Note that recurrences (\ede) and (\edg) express $\M(SC_{x,y,z}(a,b,c,m))$ in terms of five quantities of the form $\M(SC_{x',y',z'}(a,b,c,m))$, with $x'+y'+z'<x+y+z$ for each of them. By the induction hypothesis, formulas~(\eba) and (\ebb) hold for each of these five quantities. To complete the proof, it suffices to show that the expression obtained for $\M(SC_{x,y,z}(a,b,c,m))$ when substituting the formulas given by (\eba) and (\ebb) for these five quantities is equal to the right-hand side of (\eba) (for recurrence (\ede)), respectively the right-hand side of~(\ebb) (for recurrence (\edg)). Defining
$$
R_{x,y,z}(a,b,c,m):=
\cases
\text{\rm right-hand side of (\eba),}&\text{$x$, $y$, $z$ of same parity,}
\\
\\
\text{\rm right-hand side of 
(\ebb),}&\text{$x$ of parity opposite to $y$ and $z$,}
\endcases
\tag\edh
$$
this is equivalent to verifying that
$$
\spreadlines{3\jot}
\align
R_{x,y,z}(a,b,c,m)\,R_{x,y-1,z-1}(a,b,c,m)
&=
R_{y,x,z-1}(b,a,c,m)\,R_{z,y-1,x}(c,b,a,m)
\\
&
+
R_{x-1,y,z}(a,b,c,m)\,R_{x+1,y-1,z-1}(a,b,c,m),
\\
&\!\!\!\!\!\!\!\!\!\!\!\!\!\!\!\!\!\!\!\!\!\!\!\!\!\!\!\!
x,y,z\geq1,\ \text{\rm $x$ of opposite parity to $y$ and $z$}
\tag\edi
\endalign
$$
and
$$
\spreadlines{3\jot}
\align
R_{x,y,z}(a,b,c,m)\,R_{x,z-1,y-1}(a,c,b,m)
&
=
R_{z-1,x,y}(c,a,b,m)\,R_{y-1,z,x}(b,c,a,m)
\\
&
+
R_{x-1,z,y}(a,c,b,m)\,R_{x+1,y-1,z-1}(a,b,c,m),
\\
&\!\!\!\!\!\!\!\!\!\!\!\!\!\!\!\!\!\!\!\!\!\!\!\!\!\!\!\!
x,y,z\geq1,\ \text{\rm $x$, $y$ and $z$ of same parity}.
\tag\edj
\endalign
$$
We verify (\edj) first, so assume that $x$, $y$ and $z$ have the same parity. We compute the quotient of the right-hand side by the left-hand side in (\edj). Using the ``cancellation rule"
$$
\frac {\h\left(\lceil x+y\rceil\right)\,\h\left(\lfloor x+y\rfloor\right)} 
{\h\left(\lceil{x-\frac {1} {2}}\rceil+y\right)\,
\h\left(\lfloor {x-\frac {1} {2}}\rfloor+y\right)}
=\Gamma\left(\lceil x+y\rceil\right), 
\tag\edk
$$
where $x$ is an integer or a half-integer,
this quotient can be significantly simplified. Namely, after applying it numerous times, it becomes
$$
\spreadlines{3\jot}
\align
\frac{\Gamma
   (y+z+a+b+c+m) \, \Gamma
   (x+y+z+a+b+
   c+m-1) }{\Gamma
   (y+z+a+b+c+m-1) \, \Gamma
   (x+y+z+a+b+
   c+m) }\\
\times
\frac { \Gamma
   \left(\left\lceil
   \frac
   {x+y+z-1}2\right\rceil+a+b+c+m
   \right) \, \Gamma
   \left(\left\lceil
   \frac
   {x+y+z}2\right\rceil+\frac
   {a+b+
   c+m}2
   \right)} 
{  \Gamma
   \left(\left\lceil
   \frac
   {x+y+z}2\right\rceil+a+b+c+m
   \right)\, \Gamma
   \left(\left\lceil
   \frac
   {x+y+z-1}2\right\rceil+\frac
   {a+b+
   c+m}2
   \right) }
\\
+\frac{ \Gamma
   (x+a+b+c+m+1) \,\Gamma
   (x+y+z+a+b+
   c+m-1) }{ \Gamma
   (x+a+b+c+m) \, \Gamma
   (x+y+z+a+b+
   c+m) }\\
\times
\frac{\Gamma
   \left(\left\lceil
   \frac{x+1}{2}\right\rceil
   \right) \, \Gamma
   \left(\left\lceil
   \frac{x}{2}\right\rceil+\frac
   {a+b+
   c+m}2
   \right) }{\Gamma
   \left(\left\lceil
   \frac{x}{2}\right\rceil
   \right) \, \Gamma
   \left(\left\lceil
   \frac{x+1}{2}\right\rceil+\frac
   {a+b+
   c+m}2
   \right) }\\
\times
\frac{ \Gamma
   \left(\left\lceil
   \frac
   {x+y+z-1}2\right\rceil+a+b+c+m
   \right) \, \Gamma
   \left(\left\lceil
   \frac
   {x+y+z}2\right\rceil+\frac
   {a+b+
   c+m}2
   \right)}{ \Gamma
   \left(\left\lceil
   \frac
   {x+y+z}2\right\rceil+a+b+c+m
   \right)\, \Gamma
   \left(\left\lceil
   \frac
   {x+y+z-1}2\right\rceil+\frac
   {a+b+
   c+m}2
   \right) }.
\tag\edl
\endalign
$$

Let first $x,y,z$ all be even. Then the expression in (\edl) reduces to
$$
\spreadlines{3\jot}
\align
&
\frac{
   (y+z+a+b+c+m-1) }{
   (x+y+z+a+b+
   c+m-1) }
\\
&\ \ \ \ \ \ \ \ \ \ \ \ \ \ \ \ \ \ \ \ \ \ \ \ \ \ \ \ \ \ \ \ \ 
+\frac{ 
   (x+a+b+c+m)}{
   (x+y+z+a+b+
   c+m-1) }
\cdot
\frac{x /2 }{(x+a+b+c+m)/2 },
\endalign
$$
which indeed equals $1$. 

Now we assume that all of $x,y,z$ are odd.
In that case, the expression in (\edl) reduces to
$$
\spreadlines{3\jot}
\align
&
\frac{
   (y+z+a+b+c+m-1) }{
   (x+y+z+a+b+
   c+m-1) }
\cdot
\frac { 
   \left(
   \frac
   {x+y+z-1}2+\frac
   {a+b+
   c+m}2
   \right)} 
{  
   \left(
   \frac
   {x+b+c-1}2+a+b+c+m
   \right) }
\\
&\ \ \ \ \ \ \ \ \ \ \ \ \ \ \ \ \ \ \ \ \ \ \  
+\frac{ 
   (x+a+b+c+m)}{
   (x+y+z+a+b+
   c+m-1) }
\cdot
\frac{
   \left(
   \frac
   {x+y+z-1}2+\frac
   {a+b+
   c+m}2
   \right)}{ 
   \left(
   \frac
   {x+y+z-1}2+a+b+c+m
   \right) },
\endalign
$$
which also equals $1$.

\medskip
We next turn to verifying (\edi), so assume that $x$ has parity different 
from the parity of $y$ and $z$. Here, we compute the quotient of
the right-hand side by the left-hand side in (\edi). 
After applying the cancellation rule (\edk) numerous times, this quotient becomes
$$
\spreadlines{3\jot}
\align
\frac{\Gamma (y+z+a+b+c+m) \, \Gamma
    (x+y+z+a+b+c+m-1) }
{\Gamma
    (y+z+a+b+c+m-1) \,
\Gamma \left(x+y+z+a+b+c+m\right) }\\
\times
\frac{  \Gamma
    \left(\left\lceil \frac
    {x+y+z-1}2\right\rceil+a+b+c+m \right) \, \Gamma \left(\left\lceil \frac
    {x+y+z}2\right\rceil+\frac
    {a+b+c+m}2 \right)}
{  \Gamma
    \left(\left\lceil \frac
    {x+y+z}2\right\rceil+a+b+c+m \right)
 \,  \Gamma \left(\left\lceil \frac
    {x+y+z-1}2\right\rceil+\frac
    {a+b+c+m}2 \right) }\\
+
\frac {\Gamma
    (x+a+b+c+m+1) \, \Gamma (x+y+z+a+b+c+m-1) } 
{\Gamma (x+a+b+c+m) \, \Gamma
    (x+y+z+a+b+c+m) }\\
\times
\frac{
  \Gamma \left(\left\lceil
    \frac{x+1}{2}\right\rceil \right) \, 
 \Gamma \left(\left\lceil \frac{x}{2}\right\rceil+\frac
    {a+b+c+m}2
    \right) }
{  \Gamma \left(\left\lceil
    \frac{x}{2}\right\rceil \right) 
    \, \Gamma \left(\left\lceil
    \frac{x+1}{2}\right\rceil+\frac {a+b+c+m}2 \right) }\\
\times
\frac{
   \Gamma  \left(\left\lceil \frac
    {x+y+z-1}2\right\rceil+a+b+c+m \right) \, \Gamma \left(\left\lceil \frac
    {x+y+z}2\right\rceil+\frac
    {a+b+c+m}2 \right)}
{\Gamma
    \left(\left\lceil \frac
    {x+y+z}2\right\rceil+a+b+c+m \right)\,
   \Gamma \left(\left\lceil \frac
    {x+y+z-1}2\right\rceil+\frac
    {a+b+c+m}2 \right) }.
\tag\edm
\endalign
$$
Note that this happens to be identical to the expression (\edl).
Since $y+z$ is even in our case, regardless of the particular choice
of $y$ and $z$, and the expression in (\edm) depends only on $x$ and $y+z$
(and not on $y$ and $z$ separately!), we may use again the earlier
computations to deduce that this always simplifies to $1$.

This completes the verification of (\edi) and (\edj), and hence the proof. \epf

\medskip
\flushpar
{\smc Remark 1.} The special case 
$a=b=c=0$ of Theorems~{\tba} and {\tbb} was the main result of
our earlier paper \cite{\cekz} (see Theorems~1 and 2 there). The
proofs given in \cite{\cekz} used the method of factor exhaustion, and
required considerable calculations (and space!). The proof we
presented above represents in particular a great simplification 
over the proofs in~\cite{\cekz}.

\topinsert
\centerline{\mypic{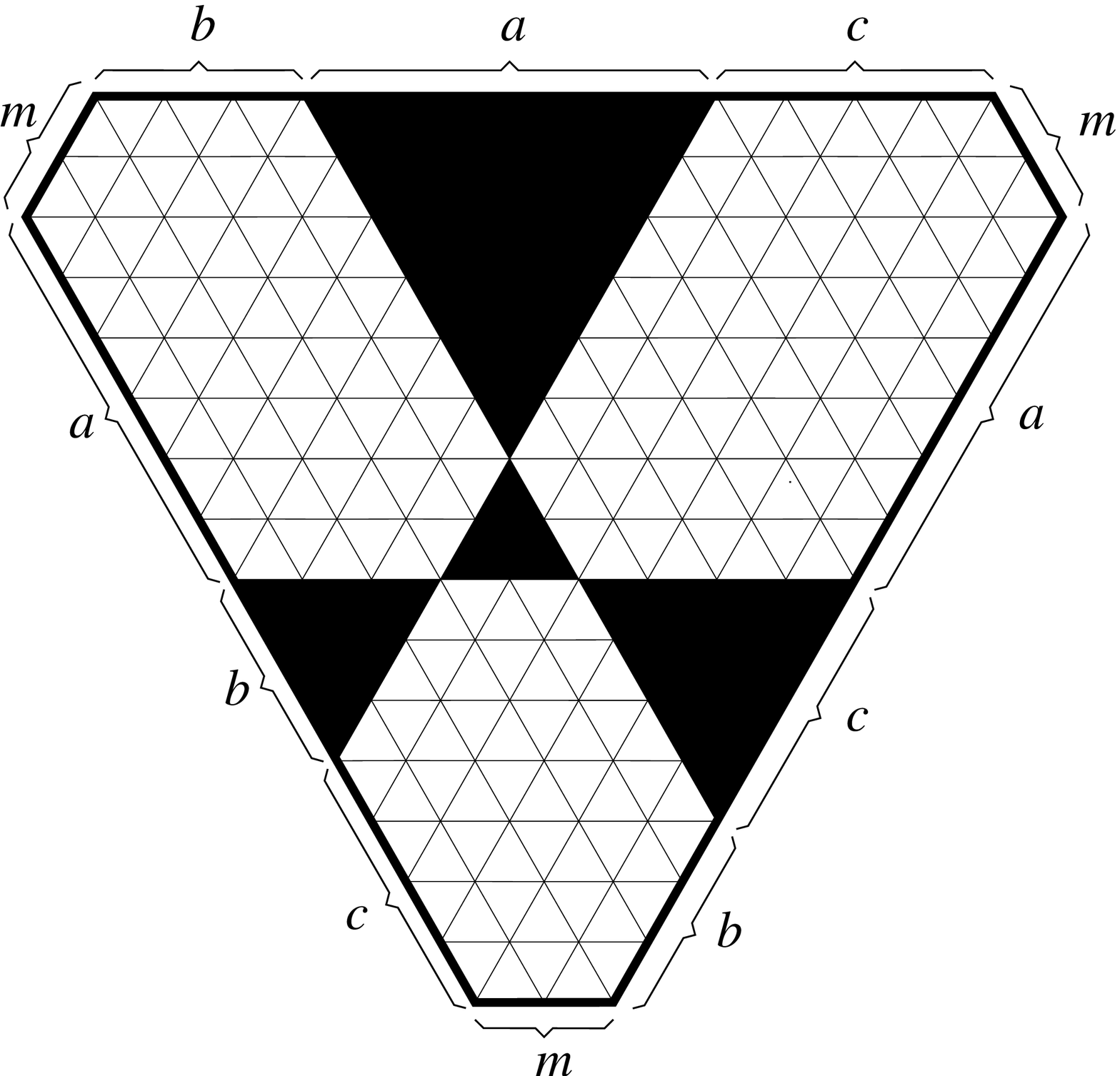}}
\medskip
\centerline{{\smc Figure~{\fdf}{. $SC_{0,0,0}(3,4,6,2)$.}}}
\endinsert

\medskip
\flushpar
{\smc Remark 2.} At the other extreme, the specialization $x=y=z=0$ of Theorems~{\tba} and~{\tbb} provides the following interesting picture (see Figure~{\fdf}). In this case, the $S$-cored hexagon $SC_{0,0,0}(a,b,c,m)$ is the disjoint union of three hexagons. It follows that
$$
\M(SC_{0,0,0}(a,b,c,m))=P(a,b,m)\,P(a,c,m)\,P(b,c,m).
\tag\edn
$$
This offers a different way of viewing Theorems~{\tba} and {\tbb} as a generalization of MacMahon's theorem~(\eaa). Namely, they can be regarded as simultaneously generalizing three applications of MacMahon's formula.

\mysec{5. Deducing Theorem {\taa}. Geometric and physical interpretation}

{\it Proof of Theorem {\taa}.}
By (\eac), the explicit definition of the expression on the left-hand side of (\ead) is
$$
\frac{\M(S^*(a,b,c,m))}{\M(S^*(a+b+c,0,0,m))}
:=
\lim_{N\to\infty}\frac{\M(SC_{N,N,N}(a,b,c,m))}{\M(SC_{N,N,N}(a+b+c,0,0,m))}.
\tag\eea
$$
In the $S$-cored hexagons on the right-hand side 
above, the $x$-, $y$- and $z$-parameters have the same parity (being all equal to $N$). Thus the number of their lozenge tilings is given by the formula provided in Theorem~{\tba}. Applying this formula one obtains, after simplifications, that
$$
\spreadlines{3\jot}
\align
&
\frac{\M(SC_{N,N,N}(a,b,c,m))}{\M(SC_{N,N,N}(a+b+c,0,0,m))}
=
\frac
{
\dfrac{\h(m)^2\h(m+a+b+c)}{\h(m+a)\,\h(m+b)\,\h(m+c)}
}
{
\dfrac{\h(a+b+c)}{\h(a)\,\h(b)\,\h(c)}
}
\\
&\ \ \ \ \ \ \ \ \ \ \ \ \ 
\times
\frac
{
\dfrac{\h(N+m)\,\h(N+m+a+b)\,\h(N+m+a+c)\,\h(N+m+b+c)}{\h(N+m+a)\,\h(N+m+b)\,\h(N+m+c)\,\h(N+m+a+b+c)}
}
{
\dfrac{\h(N)\,\h(N+a+b)\,\h(N+a+c),\h(N+b+c)}{\h(N+a)\,\h(N+b)\,\h(N+c)\,\h(N+a+b+c)}
}.
\tag\eeb
\endalign
$$
Recall that, by the Glaisher--Kinkelin formula (see \cite{\Glaish}) 
which gives the asymptotics of the Barnes $G$-function, we have
$$
\lim_{n\to\infty}
\dfrac
 {0!\,1!\,\cdots\,(n-1)!}
 {n^{\frac{n^2}{2}-\frac{1}{12}}\,(2\pi)^{\frac{n}{2}}\,e^{-\frac{3n^2}{4}}}
=
\dfrac
 {e^{\frac{1}{12}}}
 {A},
\tag\eec
$$
where $A=1.28242712...$ is a constant (called the Glaisher--Kinkelin constant). This readily implies that
$$
\lim_{N\to\infty}
\frac{\h(N)\,\h(N+a+b)\,\h(N+a+c),\h(N+b+c)}{\h(N+a)\,\h(N+b)\,\h(N+c)\,\h(N+a+b+c)}
=1.
\tag\eed
$$
As $m$ is fixed, (\eed) shows that both the numerator and the denominator of the second fraction on the right-hand side of (\eeb) approach 1 as $N\to\infty$. It follows then from (\eeb) that
$$
\lim_{N\to\infty}
\frac{\M(SC_{N,N,N}(a,b,c,m))}{\M(SC_{N,N,N}(a+b+c,0,0,m))}
=
\frac
{
\dfrac{\h(m)^2\h(m+a+b+c)}{\h(m+a)\,\h(m+b)\,\h(m+c)}
}
{
\dfrac{\h(a+b+c)}{\h(a)\,\h(b)\,\h(c)}
}.
\tag\eee
$$
The first equality in (\ead) follows from (\eea) and (\eee). The second equality in (\ead), as well as (\eae), follow then using formula (\eaa). \epf

The hexagons occurring in the symmetric form (\eae) of Theorem~{\eaa} (corresponding to the case when $m=a+b+c$) have the following geometric interpretation. Reflect each lobe of the shamrock across the vertex of the central triangle touching that lobe. This cuts out from the central triangle of the shamrock a hexagonal region, whose side-lengths are, clockwise from top, $a$, $b$, $c$, $a$, $b$, $c$. We call this the {\it central hexagon} of the shamrock. For $a=3$, $b=1$, $c=2$, this is indicated by the solid line contour in the center of Figure~{\fea}.

Another natural hexagon is obtained as follows. Consider the six pairs
of next-to-nearest edges of the central hexagon (we call two edges of a hexagon next-to-nearest if they are not incident to one another, but are both incident to another edge of the hexagon). Three of these pairs
cross at the vertices of the central triangle of the shamrock. Extend
the other three pairs as well to obtain three more crossing
points. Define the {\it medial hexagon} of the shamrock to be the
smallest lattice hexagon containing the above six crossing points. The
solid line contour on the right in Figure~{\fea} illustrates this in
the case 
$a=3$, $b=1$, $c=2$.

\topinsert
\centerline{\mypic{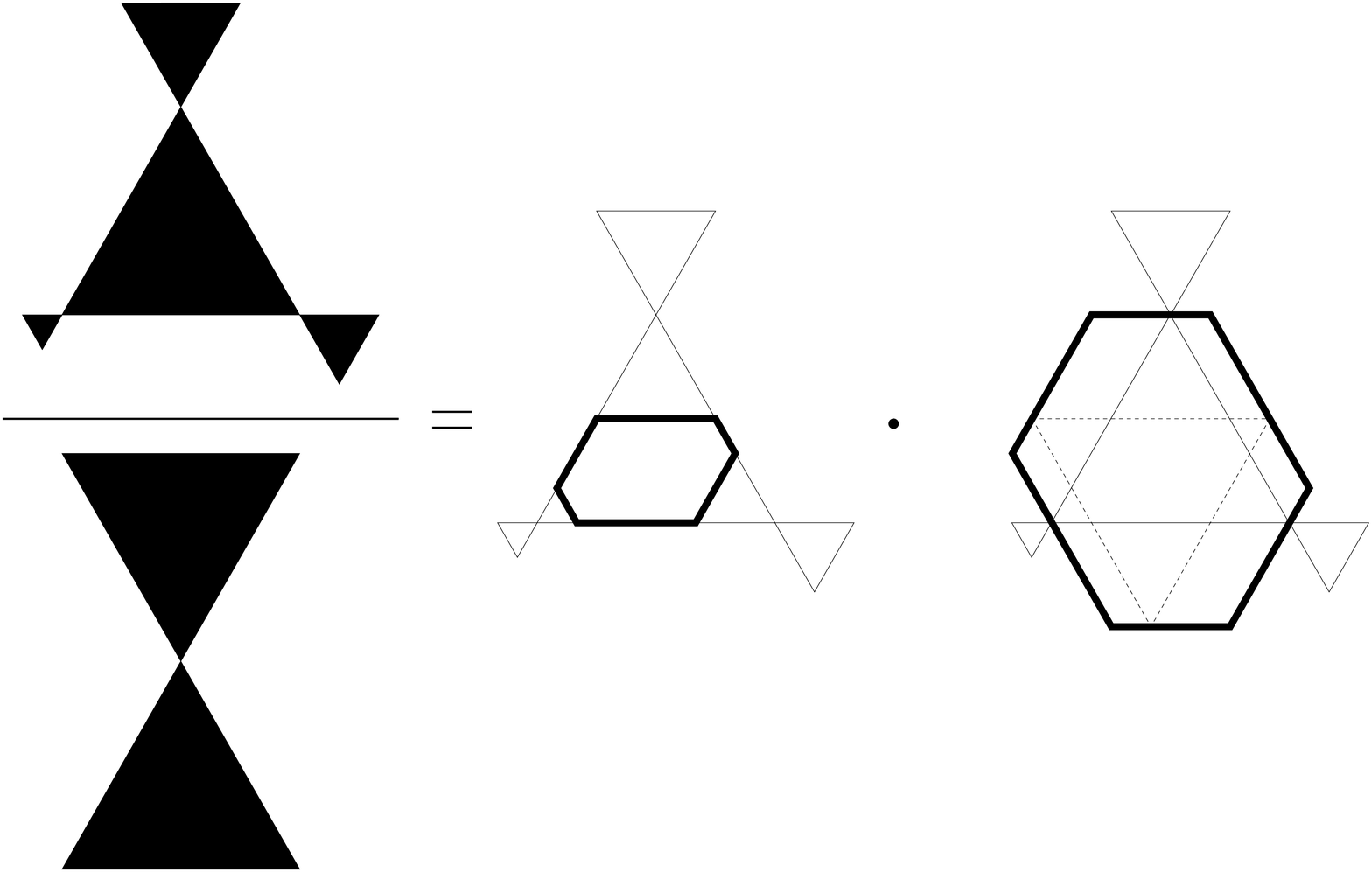}}
\medskip
\centerline{{\smc Figure~{\fea}{\rm . Geometric interpretation of Theorem~{\taa}}}} 
\centerline{\ \ \ \ \ \ \ \ \ \ \ \ \ when $m=a+b+c$; here $a=3$, $b=1$, $c=2$.}
\endinsert

Then formula~(\eae) can be expressed geometrically as shown in Figure~{\fea}: The ratio of the number of tilings of the exteriors of the two shamrocks is equal to the product of the number of tilings of the central and medial hexagons of the shamrock.

\medskip
\flushpar
{\smc Remark 3.} In the language of the first author's earlier series of papers on the correlation of holes in a sea of dimers and their connection to two dimensional electrostatics (see \cite{\sc}\cite{\ec}\cite{\ov}\cite{\ef}\cite{\gd}), the right-hand side of (\eac) is the ratio between the correlation $\omega$ defined in \cite{\sc} of the two shamrocks. Therefore (\ead) can be restated as
$$
\spreadlines{3\jot}
\align
\frac{\omega(S(a,b,c,m))}{\omega(S(a+b+c,0,0,m))}
&=
\frac{\h(a)\h(b)\h(c)\h(a+b+c+m)\h(m)^2}{\h(a+m)\h(b+m)\h(c+m)\h(a+b+c)}
\\
&=P(a,b,m)\,P(a+b,c,m),
\tag\eef
\endalign
$$
and (\eae) as
$$
\frac{\omega(S(a,b,c,a+b+c))}{\omega(S(a+b+c,0,0,a+b+c))}
=P(a,b,c)\,P(a+b,b+c,c+a).
\tag\eeg
$$
In particular, it follows that for all values of the parameters $a$,
$b$, $c$ and $m$, the ratio of the two correlations on the left-hand
side of (\eef) is an integer. There does not seem to be any simple 
a priori reason why this should be so. 

Theorem~{\taa} also has an interpretation in the light of the connection to electrostatics cited above. According to this, holes on the triangular lattice in a sea of lozenges interact (when their correlation is taken, or measured) precisely like charges in two dimensional electrostatics. Namely, the correlation of the holes is equal, in the limit of large mutual separations between the holes, to the exponential of the negative of the electrostatic energy of the two dimensional system of electrical charges obtained by regarding each hole as an electrical charge, of magnitude equal to the number of up-pointing unit triangles in the hole minus the number of down-pointing unit triangles in the hole. 

In the set-up of Theorem~{\taa}, we have four triangular holes making
up the shamrock whose correlation is considered. One readily sees that
in the above described connection to electrostatics, the lobes
correspond to charges of magnitudes $a$, $b$ and $c$, while the
central triangle corresponds to a charge of magnitude $-m$. The three
lobes should then repel one another. 
And indeed, by (\eef) the least likely configuration (i.e., the configuration of maximum energy), when one fixes $m$ and the value of the sum $a+b+c$, is the situation when one lobe contains the entire positive charge, and the other two lobes vanish.

Note that formulas~(\eef) and (\eeg) are relative correlation results, expressing the ratio of the correlation of two shamrocks. We end this remark by explaining how the correlations on the left-hand side of (\eeg) can be computed separately.

To do this, since we have formula~(\eeg), it suffices to determine $\omega(S(a+b+c,0,0,m))$. After a straightforward (if somewhat lengthy) calculation, using Stirling's formula and the Glaisher--Kinkelin formula~(\eec) to analyze the asymptotics of the resulting expressions of type~(\eaa) and (\eba), one obtains that
$$
\omega(S(a+b+c,0,0,m))
:=
\lim_{x\to\infty}
\frac{\M(SC_{x,x,x}(m,0,0,m))}{\M(\H(x+m,x+m,x+m))}
=
\frac{\sqrt{3}^{m^2}}{(2\pi)^m}
\frac{\h(m)^4}{\h(2m)}
\tag\eeh
$$
(here $\H(x+m,x+m,x+m)$ denotes the regular hexagon of side $x+m$).

\medskip
\flushpar
{\smc Remark 4.} We end this section with an interpretation of
Theorems~{\fba} and {\fbb} in terms of surfaces analogous to plane
partitions, but having a certain defect in the middle. We confine our
discussion to the case 
$m=a+b+c$ (the general case also admits an interpretation along these lines, but it is less elegant).

Recall that a lozenge tiling of a hexagon naturally lifts to a three
dimensional surface (see 
e.g.\ \cite{\DT}). This lifting works also in
the presence of a hole of charge zero (see the detailed discussion
in \cite{\ec}). For instance, the tiling shown on the left in
Figure~{\ffa} lifts to the surface shown on the right in the same
figure (note that this surface has an overhang and some precipices in
the middle). From this point of view, the case 
$m=a+b+c$ of Theorems~{\fba} and {\fbb} states that the number of such surfaces is given by the explicit product formulas~(\eba) and (\ebb).


\topinsert
\twoline{\mypic{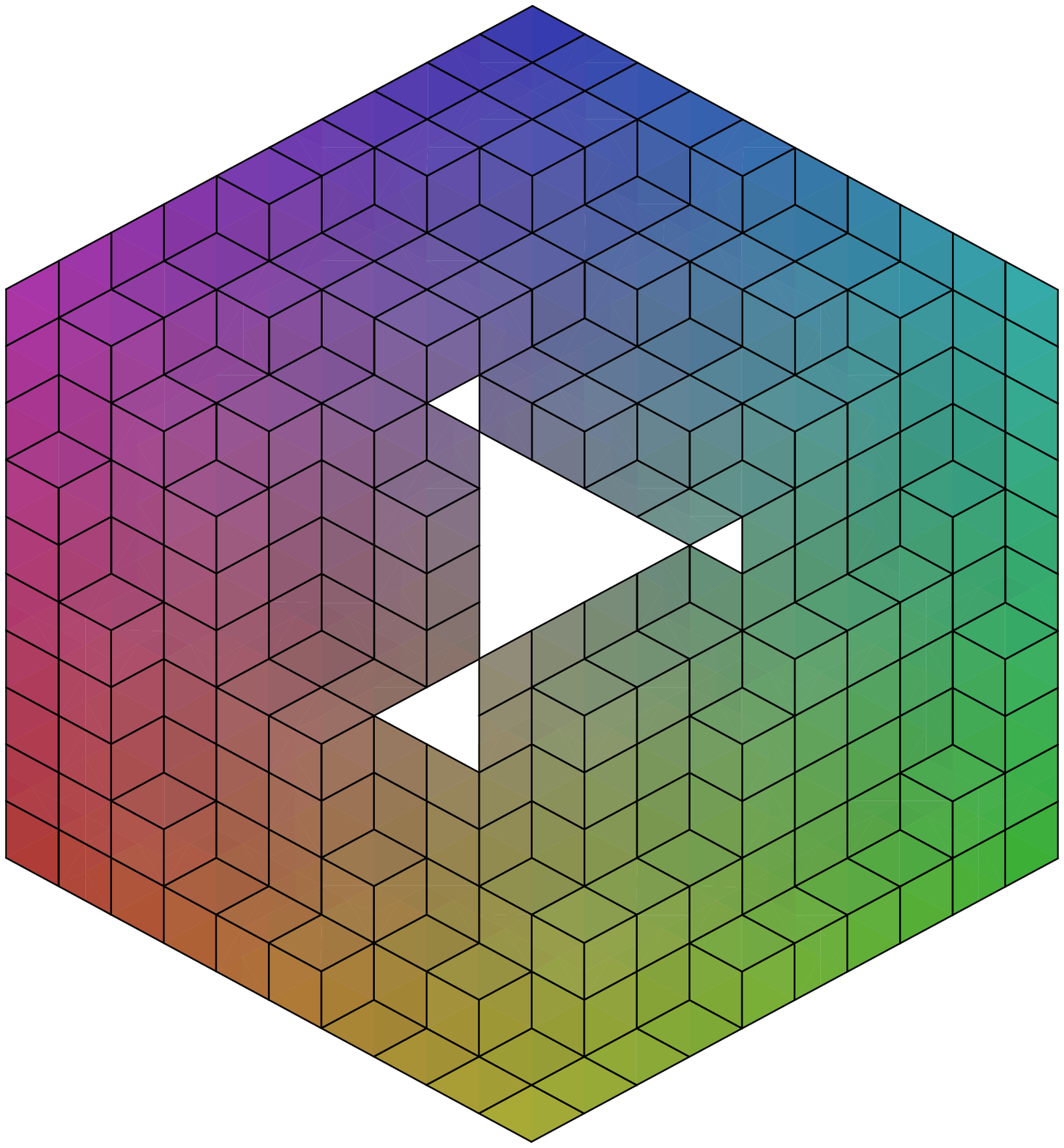}}{\mypic{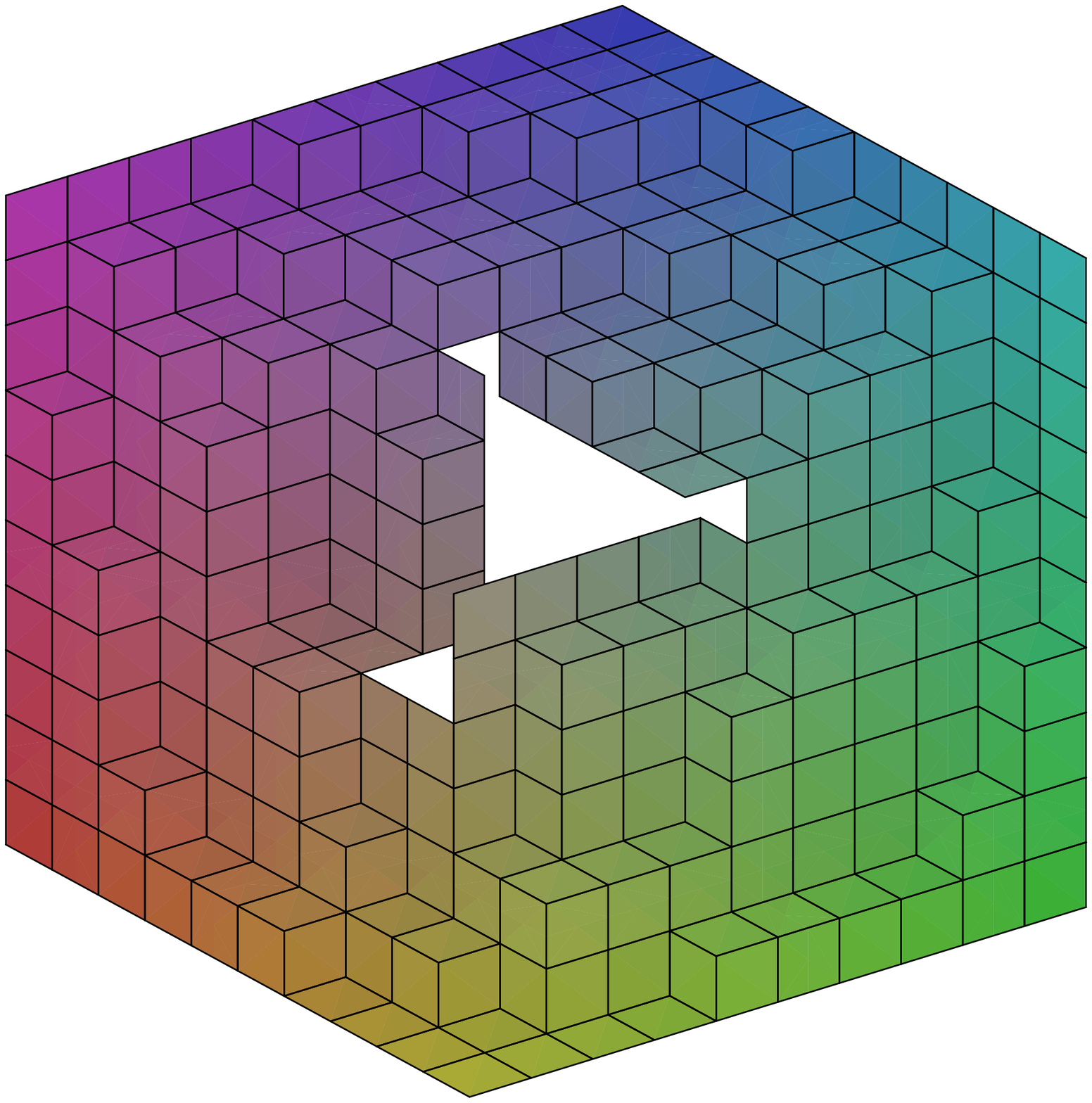}}
\centerline{Figure~{\ffa}{\rm. A tiling of $S_{6,6,6}(1,2,1,4)$ and its lifting.}}
\endinsert







\mysec{6. Concluding remarks}

We have seen in this paper the great power of Kuo's graphical
condensation method, when applied to situations in which explicit
conjectured formulas can be found (for another application of the same
method, see \cite{\anglepap}). The proofs presented here represent a
great simplification of the proofs of Theorems~1 and 2
of \cite{\cekz}, which are the special case 
$a=b=c=0$ of Theorems~{\tba} and {\tbb} of this paper. Furthermore, the same method can be used to tackle Conjectures~1 and 2 of \cite{\cekz}, which seem otherwise quite redoubtable, and also to extend them to $S$-cored hexagons (details will appear in a separate paper).






\mysec{References}
{\openup 1\jot \frenchspacing\raggedbottom
\roster

\myref{\And}
  G. E. Andrews, Plane partitions, III: The weak Macdonald
conjecture, {\it Invent. Math.} {\bf 53} (1979), 193--225.

\myref{\Andtwo}
  G.E. Andrews, Plane partitions, V: The T.S.S.C.P.P. conjecture, {\it J. Combin. Theory Ser.~A} {\bf 66} (1994), 28–-39.


\myref{\cekz}
  M. Ciucu, 
T. Eisenk\"olbl, C. Krattenthaler and D. Zare,
Enumeration of lozenge tilings of hexagons with a central triangular hole,
{\it J. Combin. Theory Ser.~A} {\bf 95} (2001), 251--334.

\myref{\sc}
  M. Ciucu, A random tiling model for two dimensional electrostatics, 
{\it Mem. Amer. Math. Soc.} {\bf 178} (2005), no. 839, 1--106.

\myref{\ec}
  M. Ciucu, The scaling limit of the correlation of holes on the triangular lattice
with periodic boundary conditions, {\it Mem. Amer. Math. Soc.} {\bf 199} (2009),
no. 935, 1-100.

\myref{\ov}
  M. Ciucu, Dimer packings with gaps and electrostatics, {\it Proc. Natl. Acad. Sci.
USA} {\bf 105} (2008), 2766-2772.

\myref{\ef}
  M. Ciucu, The emergence of the electrostatic field as a Feynman sum in random
tilings with holes, {\it Trans. Amer. Math. Soc.} {\bf 362} (2010), 4921-4954.

\myref{\gd}
  M. Ciucu, The interaction of collinear gaps of arbitrary charge in a two
dimensional dimer system, 2012, 
{\tt ar$\chi$iv:1202.1188}.

\myref{\anglepap}
  M. Ciucu and I. Fischer, A triangular gap of size two in a sea of dimers on a
$60^\circ$ angle, preprint, July 2012, submitted.

\myref{\CLP}
  H. Cohn, M. Larsen, and J. Propp, The shape of a typical boxed plane partition,
{\it New York J. of Math.} {\bf 4} (1998), 137--165.

\myref{\DT}
  G. David and C. Tomei, The problem of the calissons,
{\it Amer\. Math\. Monthly} {\bf 96} (1989), 429--431.


\myref{\Glaish}
  J. W. L. Glaisher, On certain numerical products in which the
  exponents depend upon the numbers, {\it  Messenger Math.} {\bf 23}
  (1893), 145--175. 

\myref{\KaKZAC}
C. Koutschan, M. Kauers and D. Zeilberger, A proof of George
Andrews' and David Robbins' $q$-TSPP-conjecture, {\it Proc.\
Natl.\ Acad.\ Sci.\ USA} {\bf 108} (2011), 2196--2199.

\myref{\Kuo}
E. H. Kuo,
Applications of graphical condensation for enumerating matchings and tilings. Theoret. Comput. Sci. 319 (2004), no. 1-3, 29–57.

\myref{\Kup}
  G. Kuperberg, Symmetries of plane partitions and the permanent-de\-ter\-mi\-nant
method, {\it J. Combin. Theory Ser.~A} {\bf 68} (1994), 115--151.

\myref{\MacM}
  P. A. MacMahon, Memoir on the theory of the partition of numbers---Part V. 
Partitions in two-dimensional space, {\it Phil. Trans. R. S.}, 1911, A.

\myref{\Sta}
  R. P. Stanley, Symmetries of plane partitions, {\it J. Combin. Theory. Ser.~A}
{\bf 43} (1986), 103--113.

\myref{\Ste}
  J. R. Stembridge, The enumeration of totally symmetric plane partitions,
{\it Adv. in Math.} {\bf 111} (1995), 227--243.

\endroster\par}

\enddocument